\documentclass[12pt,leqno]{amsart}

\usepackage{amsmath,amscd,amsthm,amsxtra,amssymb}
\usepackage{epsfig,graphics,color,colortbl}
\usepackage{amssymb,latexsym}
\usepackage{mathrsfs}
\usepackage[poly,all]{xy}
\usepackage{marginnote}
\usepackage{xspace}
\usepackage[colorlinks=true, pdfstartview=FitV, linkcolor=blue,citecolor=blue,urlcolor=blue]{hyperref}

\usepackage{yfonts}
\usepackage{enumerate}

\usepackage[usenames,dvipsnames,svgnames,table]{xcolor}

\usepackage[normalem]{ulem}  

\usepackage[colorinlistoftodos]{todonotes}

\newdir{ >}{{}*!/-10pt/@{>}}

\allowdisplaybreaks[3]

\newlength{\mylength}
\renewcommand{\le}{\leqslant}
\renewcommand{\ge}{\geqslant}

\theoremstyle{plain}
\newtheorem{thm}{Theorem}[section]
\newtheorem{df}[thm]{Definition}
\newtheorem{prop}[thm]{Proposition}
\newtheorem{coro}[thm]{Corollary}
\newtheorem{lem}[thm]{Lemma}
\newtheorem{conj}[thm]{Conjecture}

\theoremstyle{definition}
\newtheorem{ex}[thm]{Example}
\newtheorem{remark}[thm]{Remark}
\newtheorem{definition}[thm]{Definition}
\newtheorem{question}[thm]{Qusetion}

\newcommand{\nc}{\newcommand}

\newenvironment{answer}
{\noindent{\bf Answer}\hs{1ex}}
{\hfill \qedsymbol}
\nc{\Prop}{\begin{prop}}
\nc{\enprop}{\end{prop}}
\nc{\Lemma}{\begin{lem}}
\nc{\enlemma}{\end{lem}}
\nc{\Exam}{\begin{ex}}
\nc{\enexam}{\end{ex}}
\nc{\Th}{\begin{thm}}
\nc{\enth}{\end{thm}}
\nc{\Def}{\begin{definition}}
\nc{\edf}{\end{definition}}
\nc{\Conj}{\begin{conj}}
\nc{\enconj}{\end{conj}}
\nc{\Quest}{\begin{question}}
\nc{\enquest}{\end{question}}
\nc{\Rem}{\begin{remark}}
\nc{\enrem}{\end{remark}}
\nc{\Ans}{\begin{answer}}
\nc{\enans}{\end{answer}}

\newenvironment{red}
{\relax\color{red}}
{\hspace*{.5ex}\relax}

\newcommand{\ber}{\begin{red}}
\newcommand{\er}{\end{red}}

\nc{\berm}{\ber{}\marginnote{\fbox{\scshape\lowercase{M}}}}
\nc{\berMH}{\ber{}\marginnote{\fbox{\scshape\lowercase{MH}}}}
\nc{\berE}{\ber{}\marginnote{\fbox{\scshape\lowercase{E}}}}
\newenvironment{blue}
{\relax\color{Dandelion}}
{\hspace*{.5ex}\relax}

\newcommand{\beb}{\begin{blue}}
\newcommand{\eb}{\end{blue}}

\nc{\on}{\operatorname}

\newcommand{\C}{{\mathbb C}}
\newcommand{\Q}{\mathbb {Q}}
\newcommand{\Z}{\ms{2mu}{\mathbb Z}}

\newcommand{\D}{\mathscr{D}\ms{1mu}}

\newcommand{\one}{{\bf{1}}}
\newcommand{\seteq}{\mathbin{:=}}

\newcommand{\To}[1][{\hs{0.8ex}}]{\xrightarrow{\ms{7mu}{#1}\ms{7mu}}}

\newcommand{\g}{\ms{1mu}\mathfrak{g}\ms{1mu}}
\newcommand{\n}{\mathfrak{n}}

\newcommand{\Hom}{\operatorname{Hom}}
\newcommand{\HOM}{\on{\mathrm{H{\scriptstyle OM}}}}

\newcommand{\isoto}[1][]{\mathop{\xrightarrow%
[{\raisebox{.3ex}[0ex][.3ex]{$\scriptstyle{#1}$}}]%
{{\raisebox{-.6ex}[0ex][-.6ex]{$\mspace{2mu}\sim\mspace{2mu}$}}}}}

\newcommand{\Mod}{\on{Mod}}
\newcommand{\gmod}{\text{-}\mathrm{gmod}}

\newcommand{\F}{\mathscr{F}}

\newcommand{\conv}[1][]{
\underset{\raisebox{.5ex}{$\scriptstyle{#1}$}}{\mathbin{\scalebox{1.1}{$\mspace{1.5mu}\circ\mspace{1.5mu}$}}}}
\newcommand{\hconv}{\mathbin{\scalebox{.9}{$\nabla$}}}

\newcommand{\sconv}{\mathbin{\scalebox{.9}{$\Delta$}}}

\renewcommand{\Im}{\on{Im}}
\newcommand{\de}{\on{\textfrak{d}}}

\newcommand{\cmA}{\cartan}  
\newcommand{\wlP}{\mathsf{P}}   
\newcommand{\rlQ}{\mathsf{Q}}   
\newcommand{\weyl}{\mathsf{W}}  
\newcommand{\prD}{\Delta_+}            
\newcommand{\sg}{\mathfrak{S}}   
\newcommand{\Po}{\wlP}

\nc{\prt}{\prD}
\nc{\qQ}{Q}
\newcommand{\bQ}{\overline{\qQ}}

\newcommand\Aq[1][{\g^+}]{A_q(#1)}

\newcommand{\wt}{\mathrm{wt}} 		
\newcommand{\bR}{\mathbf{k}} 		
\nc{\corp}{\bR}
\newcommand{\catC}{ \mathscr{C}}  	
\newcommand{\tcatC}{ \widetilde{\mathscr{C}}}  	
\newcommand{\catT}{ \mathcal{T}}  	
\newcommand{\lT}{ \widetilde{\mathcal{T}}}  	

\newcommand{\Qr}{\mathrm{Q}^{\ms{1.5mu}\mathrm{r}}} 
\newcommand{\Ql}{\mathrm{Q}^{\ms{1.5mu}\mathrm{l}}} 
\newcommand{\dM}{ \mathsf{M }}              
\newcommand{\Mw}{\ms{1mu}\dM(w\La,\La)}      
\newcommand{\Mv}{ \dM(v\La,\La)}      
\newcommand{\Mwv}{ \dM(w\La,v\La)}      
\newcommand{\dC}{ \mathsf{C }}              
\newcommand{\mD}{ \mathsf{D}}              
\newcommand{\gW}{\mathsf{W}}
\newcommand{\sgW}{\mathsf{W}^*}
\newcommand{\ep}{\varepsilon}  		

\newcommand{\Ht}{\mathrm{ht}} 		
\newcommand{\nR}{\mathrm{R}^{\mathrm{norm}}} 		
\newcommand{\RR}{\mathrm{R}} 				
\newcommand{\coR}{\mathrm{R}^{\mathrm{l}}}				
\newcommand{\coRr}{\mathrm{R}^{\mathrm{r}}} 	
\newcommand{\coRl}{\mathrm{R}^{\mathrm{l}}} 	
\newcommand{\La}{\Lambda} 			
\newcommand{\tLa}{\widetilde{\Lambda}} 			
\newcommand{\Dd}{\text{ \textfrak{d}}} 			
\newcommand{\Res}{\mathrm{Res}\ms{1mu}} 			

\nc{\Ma}{{\ms{1.5mu}\mathsf{M}}}
\nc{\Na}{\mathsf{N}}
\nc{\Xa}{\mathsf{X}}
\nc{\Ya}{\mathsf{Y}}
\nc{\Laa}{\mathsf{L}}

\newcommand{\z}[1][{\Ma}]{{z_{#1}}}

\newcommand{\gr}{\mathrm{gr}}
\newcommand{\triv}{{\mathbf{1}}}   				
\newcommand{\id}{\ms{2mu}{\mathsf{id}}\ms{1mu}}   				
\newcommand{\Ds}{{\mathcal{D}}}   				
\newcommand{\Hm}{{\mathrm{H}}}   				
\newcommand{\gHm}{\mathrm{H}^{\gr}}   				
\newcommand{\dphi}{{\phi}}   				
\newcommand{\gH}{\mathrm{H}}   				
\newcommand{\gL}{\mathrm{L}}   				
\newcommand{\gzeta}{\zeta^\gr}   				
\newcommand{\gPsi}{\Psi^\gr}   				
\newcommand{\gtPsi}{\widetilde{\Psi}^\gr}   				
\newcommand{\gT}{T^\gr}   				
\newcommand{\gtT}{\widetilde{T}^\gr}   				

\newcommand{\lG}{\Gamma}   					

\nc{\be}{\begin{enumerate}}
\newcommand{\bnum}{\be[{\rm(i)}]}
\newcommand{\bna}{\be[{\rm(a)}]}

\newcommand{\rtl}{\rlQ}

\newcommand{\etens}{\boxtimes}

\newcommand{\rmat}[1]{\ms{1mu}{\mathbf{r}}_%
{\mspace{-2mu}\raisebox{-.6ex}{${\scriptstyle{#1}}$}}}

\newcommand{\shc}{{\ms{2mu}\mathcal{C}}}

\newcommand{\tC}{\widetilde{\catC}}

\newcommand{\Ob}{\on{Ob}}

\nc{\ms}{\mspace}
\nc{\cl}{\colon}
\nc{\ro}{{\rm (}}
\nc{\rf}{{\rm )}\xspace}
\nc{\noi}{\noindent}
\nc{\bl}{\bigl(}
\nc{\br}{\bigr)}

\newenvironment{myequationn}
{\relax\setlength{\arraycolsep}{1pt}\begin{eqnarray*}}
{\end{eqnarray*}}

\newenvironment{myequation}
{\relax\setlength{\arraycolsep}{1pt}\begin{eqnarray}}
{\end{eqnarray}}

\nc{\eq}{\begin{myequation}}
\nc{\eneq}{\end{myequation}}
\nc{\eqn}{\begin{myequationn}}
\nc{\eneqn}{\end{myequationn}}

\newenvironment{myarray}[1]{\relax\setlength{\arraycolsep}{1pt}
\begin{array}{#1}}{\end{array}\relax}

\newcommand{\ba}{\begin{myarray}}
\newcommand{\ea}{\end{myarray}}

\nc{\hs}{\hspace*}
\nc{\vs}{\vspace*}
\nc{\set}[2]{\left\{{#1}\mid{#2}\right\}}
\nc{\snoi}{\smallskip\noi}
\nc{\mnoi}{\medskip\noi}
\nc{\al}{\alpha}
\nc{\rmz}{\setminus\{0\}}
\nc{\tens}
[1][]{\mathbin{\otimes{\ms{2mu}}_{\raise1.5ex\hbox to-.1em{}#1}}}
\nc{\vphi}{\varphi}
\nc{\ee}{\end{enumerate}}
\nc{\la}{\lambda}
\nc{\bc}{\begin{cases}}
\nc{\ec}{\end{cases}}
\nc{\qtq}[1][and]{\quad\text{#1}\quad}
\nc{\qt}[1]{\quad\text{#1}}
\nc{\dual}{{\displaystyle{\ms{1mu}\star}}}
\nc{\wle}{\preceq}
\nc{\epito}{\twoheadrightarrow}
\nc{\epiTo}[1][]{\xymatrix@C=4ex{{}\ar@{->>}[r]^-{#1}&{}}}
\nc{\Proof}{\begin{proof}}
\nc{\lan}{\langle}
\nc{\ran}{\rangle}
\nc{\ang}[1]{\lan{#1}\ran}
\nc{\QED}{\end{proof}}
\nc{\soplus}{\scalebox{.65}{\raisebox{.2ex}{$\displaystyle\bigoplus$}}}
\nc{\eps}{\varepsilon}
\nc{\supp}{\on{supp}}
\nc{\sct}{strongly commute\xspace}
\nc{\scts}{strongly commutes\xspace}
\nc{\bce}{\eta}			
\nc{\height}[1]{\on{ht}(\ms{.5mu}{#1}\ms{.5mu})}
\nc{\braid}{{\ms{1mu}\mathrm{br}}}
\nc{\gp}{\mathfrak{p}}
\nc{\wtl}{\wlP}
\nc{\ra}{real and admits an affinization}
\nc{\ras}{real and admit affinizations}
\nc{\Cor}{\begin{coro}}
\nc{\encor}{\end{coro}}
\nc{\shf}{\mathcal{F}}
\nc{\Cw}[1][{w}]{\catC_{{#1}}}
\nc{\tCw}[1][{w}]{\widetilde{\catC}_{{#1}}} 
\nc{\tCwv}[1][{w,v}]{\widetilde{\catC}_{{#1}}} 
\nc{\akew}[1][1ex]{\rule[-1ex]{#1}{0ex}}
\nc{\ake}[1][2ex]{\rule[-1ex]{0ex}{#1}}
\nc{\akete}[1][-1ex]{\rule[{#1}]{0ex}{1ex}}
\nc{\tRm}{(R\gmod)\widetilde{\mbox{$\ake[2.5ex]\akew[.9ex]$}}}
\nc{\monoTo}[1][]{\xymatrix{\ar@{>->}[r]^-{{#1}}&}}
\nc{\monoto}[1][]{\rightarrowtail}
\nc{\tX}{\widetilde{X}}
\nc{\corps}{\corp}
\nc{\tL}{\widetilde{L}}
\nc{\prtl}{\rtl_+}
\nc{\nrtl}{\rtl_-}
\nc{\tK}{\widetilde{K}}
\nc{\tep}{\widetilde\ep}
\nc{\teps}{\widetilde\ep}
\nc{\tmu}{\widetilde \mu} 
\nc{\teta}{\widetilde\eta}
\nc{\ga}{\mathfrak{a}}
\nc{\scbul}{{\,\raise1pt\hbox{$\scriptscriptstyle\bullet$}\,}}
\nc{\bwr}{\mbox{\large$\wr$}}
\nc{\tR}{{\widetilde{\mathrm{R}}}}
\nc{\lS}{\mathsf{S}}
\nc{\lZ}{\mathcal{Z}}
\nc{\prolim}[1][]{\mathop{\varprojlim}\limits_{{#1}}}
\nc{\sym}{\sg}
\newcounter{myc}

\nc{\txi}{\tilde{\xi}}
\nc{\rl}{\rlQ}
\nc{\sfC}{\mathsf{C}}
\nc{\cor}{{\ms{1mu}\mathbf{k}\ms{1mu}}}
\nc{\Pro}{\on{Pro}}

\nc{\hM}{\widehat{\mathsf{M}}}
\nc{\aff}{\mathrm{aff}}
\nc{\rDa}{{\mathscr{D}_\aff}}
\nc{\st}[1]{\left\{{#1}\right\}}
\nc{\W}{\mathsf{W}}
\nc{\rt}{\Delta}
\nc{\pwtl}{\wtl_+}
\nc{\rev}{{\mathrm{rev}}}
\nc{\E}[1]{\mathrm{E}_{{#1}}\ms{1mu}}
\nc{\Es}[1]{\mathrm{E}^*_{{#1}}\ms{1mu}}
\nc{\Qt}{\mathscr{Q}}
\nc{\Ctr}{\mathsf{C}}
\nc{\Ctrs}{{\mathsf{C}^*}}
\nc{\Dynkin}{\Delta}
\nc{\cartan}{\mathsf{C}}
\nc{\sfc}{\mathsf{c}}
\nc{\sfa}{\mathsf{a}}
\nc{\SW}{\mathrm{K}}
\nc{\hSW}{\widehat{\mathrm{K}}}
\nc{\refl}{\mathscr{S}}
\nc{\Rre}{\mathrm{R}^{\mathrm{ren}}}
\nc{\Rpre}{\mathrm{R}^{\mathrm{ren}\;'}}
\nc{\bRre}{\ol{\mathrm{R}}^{\ms{2mu}\mathrm{ren}}}
\nc{\sfd}{\ms{1mu}\mathsf{d}\ms{1mu}}
\nc{\shm}{\mathcal{M}}
\nc{\sht}{\mathcal{T}}
\nc{\rank}{\mathrm{rank}}
\nc{\Da}{{\D}\ms{-2.8mu}\raisebox{-.35ex}{$\scriptstyle\mathrm{aff}$}}
\nc{\lDa}{\Da^{-1}}
\nc{\bchi}{{\scalebox{.9}{\mbox{$\mathscr{E}$}}}}
\nc{\bchis}{\bchi{}^{\ms{2mu}*}}
\nc{\Daf}{\mathscr{D}}
\nc{\Laf}{\mathscr{L}}
\nc{\tLaf}{\widetilde{\Laf}}
\nc{\wtaf}{\mathscr{W}\ms{-3mu}{\mathit{t}}}
\nc{\res}[1][]{\mathop\star\limits_{\raisebox{.4ex}{$\scriptstyle #1$}}\ms{2mu}}
\nc{\hchi}{\widehat{\chi}}
\nc{\convaff}{\mathop{\scalebox{1.1}{$\mspace{1.5mu}\circ\mspace{1.5mu}$}}\limits}
\nc{\cvb}{CVB\xspace}
\nc{\svelt}{essentailly samll\xspace}
\nc{\Proc}{\on{Pro}_{\mathrm{coh}}}
\nc{\sha}{\mathcal{A}}
\nc{\Ker}{\on{Ker}}
\nc{\Coker}{\on{Coker}}
\nc{\Aff}[1][z]{\on{Aff}_{\ms{1mu}#1}}
\nc{\scb}{\scalebox}
\nc{\afr}{affreal\xspace}
\nc{\epifrom}{\ms{-5mu}\xymatrix@C=3ex{{}&{}\ar@{->>}[l]}\ms{-5mu}}
\nc{\Mid}{\bigm|}
\nc{\ol}{\overline}
\nc{\bpsi}{\ol{\psi}}
\nc{\Rat}[1][z]{\on{Raff}_{\ms{1mu}#1}}
\newcommand{\indlim}[1][]{\mathop{\varinjlim}\limits_{#1}}
\nc{\inddlim}{\mathop{\mbox{``{$\ms{1mu}\varinjlim$}''}}\limits}
\nc{\Rmat}{\mathrm{R}\ms{1mu}}
\nc{\Runi}{\mathrm{R}^{\mathrm{univ}}}
\nc{\Modg}{\mathrm{Modg}}
\nc{\KO}{quasi-rigid\xspace}
\nc{\hF}{\widehat{\F}}
\nc{\Modc}{\Mod_{\mathrm{coh}}}
\nc{\e}{\mathrm{e}}
\nc{\Idx}{\mathsf{\Lambda}}
\nc{\hA}{\widehat{A}}
\nc{\prood}{\mathop{\text{``}\prod\text{''}}\limits}

\nc{\hrefl}{\widehat{\mathscr{S}}}
\nc{\ev}{\mathrm{ev}}
\nc{\coev}{\mathrm{coev}}
\nc{\ihom}{\mathcal{H}om}
\nc{\tY}{\widetilde{Y}}
\nc{\tensz}{\tens[z]\ms{-3.5mu}}

\nc{\tRre}{\widetilde{\mathrm{R}}^{\mathrm{ren}}}
\nc{\dg}{\mathbf{\lambda}}
\nc{\htens}{\hconv}
\nc{\stens}{\sconv}
\nc{\Modgc}{\mathrm{Modg}_{\mathrm{coh}}}
\nc{\Aut}{\mathrm{Aut}}
\nc{\Rd}[1][\dg]{R_{#1}\gmod}
\nc{\nn}{\nonumber}
\nc{\Dual}{\mathrm{D}\ms{1mu}}
\nc{\DA}[1][A]{\ms{1mu}\mathrm{D}_{{#1}}}
\nc{\DmA}[1][A]{\ms{1mu}\Dual^{-1}_{{#1}}}
\nc{\tensa}{\tens[A]\ms{-3mu}}
\nc{\tensc}{\tens[{\ms{3mu}\cor}]\ms{-3mu}}

\nc{\ble}{\preccurlyeq}
\nc{\bge}{\succcurlyeq}
\nc{\Cwv}{\catC_{w,v}}

\nc{\afn}{affine object\xspace}
\nc{\afns}{affine objects\xspace}
\nc{\subafn}{affine subobject\xspace}
\nc{\Afns}{Affine objects\xspace}
\nc{\mr}{\mathrm{r}}
\nc{\ml}{\mathrm{l}}
\nc{\LQ}{\mathscr{L}}
\nc{\RQ}{\mathscr{R}}
\nc{\tM}{\tilde{M}}
\nc{\tN}{\tilde{N}}
\nc{\convz}[1][z]{\underset{#1}{\circ}}

\numberwithin{equation}{section}
\title{Localizations for quiver Hecke algebras III}

\author[M. Kashiwara]{Masaki Kashiwara}
\thanks{The research of M.\ Kashiwara
was supported by Grant-in-Aid for Scientific Research (B) 20H01795,
Japan Society for the Promotion of Science.}
\address[M. Kashiwara]{ 
Kyoto University Institute for Advanced Study, Research Institute
for Mathematical Sciences, Kyoto University, Kyoto 606-8502, Japan
\& Korea Institute for Advanced Study, Seoul 02455, Korea }
\email{masaki@kurims.kyoto-u.ac.jp}

\author[M. Kim]{Myungho Kim}
\address[M. Kim]{Department of Mathematics, Kyung Hee University, Seoul 02447,   Korea}
\email{mkim@khu.ac.kr}
\thanks{The research of M.\ Kim was supported by the National Research Foundation of
Korea (NRF) Grant funded by the Korea government(MSIP)
(NRF-NRF-2022R1F1A1076214 and NRF-2020R1A5A1016126).}

\author[S.-j. Oh]{Se-jin Oh}
\thanks{ The research of S.-j.\ Oh was supported by the Ministry of Education of the Republic of Korea and the National Research Foundation of Korea (NRF-2022R1A2C1004045).}
\address[S.-j. Oh]{Department of Mathematics,  Sungkyunkwan University,  Suwon 16419,  Korea}
\email{sejin092@gmail.com}

\author[E. Park]{Euiyong Park}
\thanks{The research of E.\ Park was supported by the National Research Foundation of Korea (NRF) Grant funded by the Korea Government(MSIP)(NRF-2020R1A5A1016126).}
\address[E. Park]{Department of Mathematics, University of Seoul, Seoul 02504,   Korea}
\email{epark@uos.ac.kr}


\makeatletter
\@namedef{subjclassname@2020}{\textup{2020} Mathematics Subject Classification}
\makeatother


	\date{August 18, 2023}

\begin{document}

\maketitle

\begin{abstract}
	
Let $R$ be a quiver Hecke algebra, and let $\catC_{w,v}$ be the category of finite-dimensional graded $R$-module categorifying a $q$-deformation of the doubly-invariant algebra $^{N'(w)} \C[N] ^{N(v)} $. In this paper,
we prove that
the localization $\tcatC_{w,v}$ of the category $\catC_{w,v}$ can be
obtained as the localization by right braiders arising from  determinantial modules. 
As its application,  we show several interesting properties of the localized category $\tcatC_{w,v} $ including the right rigidity.
\end{abstract}

\tableofcontents

\section{Introduction}

This is the third of our series of papers on \emph{localizations} for \emph{quiver Hecke algebras} (\cite{KKOP21, KKOP22}). 
Let $R$ be a quiver Hecke algebra associated with a quantum group $U_q(\g)$  and let $R\gmod$ be the category of finite-dimensional graded $R$-modules.  
The monoidal category $R\gmod$ has a rich and interesting structure,  and there are various successful applications including \emph{monoidal categorification} for quantum cluster algebras (\cite{KKKO18}).  
For an element $w$ of the Weyl group $\weyl$,  there is an interesting subcategory $\catC_w$ of $R\gmod$ defined by using $w$ (see Section \ref{Sec: Cwv}  for the precise definition). 
It was proved  that the Grothendieck ring $ K( \catC_w)$ of $\catC_w$ is isomorphic to the \emph{quantum unipotent coordinate ring} $A_q(\n(w))$ associated with $w$ (see \cite{KKKO18} and see also \cite{KKOP18}).  The algebra $A_q(\n(w))$  can be understood as a $q$-deformation of the coordinate ring of the unipotent subgroup $N(w)$, and  $A_q(\n(w))$ has a quantum cluster algebra structure. It was proved in \cite{KKKO18} that the category $\catC_w$ gives a monoidal categorification of $A_q(\n(w))$ as a quantum cluster algebra when $\g$ is symmetric.

In the previous works (\cite{KKOP21, KKOP22}) by the authors,  a localization procedure using \emph{left braiders} was developed and
applied to the subcategories $\catC_w$.
For a  graded monoidal category $\catT =\soplus_{\la \in\La} \catT_\la$, 
a \emph{graded left braider} in $\catT$  is  a triple $(C, \coRl_C, \phi)$ such that $C$ is an object, $\phi\cl \La \to \Z$ is a $\Z$-linear map, and  
\begin{align} \label{Eq: lb}
\coRl_{C} (X) \cl  C \tens X  \to q^{-\phi(\la)}\tens  X \tens C
\end{align}
is a morphism functorial in $X \in \catT_\la$ such that they satisfy certain compatibility (see \cite[Section 2.3]{KKOP21}).
In \cite{KKOP21}, we proved that there exist the graded left braiders in $R\gmod$
$$
( \dM(w\La_i,\La_i),  \coR_{\dM(w\La_i,\La_i)}, \dphi^l_{\dM(w\La_i,\La_i)} ) \qquad \text{ for any $i\in I$ }
$$ 
using the \emph{determinantial modules} $\dM( w\La_i, \La_i )$, the morphism induced from \emph{R-matrices}: for any $X \in R\gmod$,
$$
\coR_{\dM(w\La_i,\La_i)}(X) \cl \dM(w\La_i,\La_i) \conv  X  \longrightarrow q^{- \dphi^l_{\dM(w\La_i,\La_i)} ( \wt(X) ) }  X  \conv   \dM(w\La_i,\La_i),  
$$
and  $\dphi^l_{\dM(w\La_i,\La_i)} ( \beta) = - (w \La_i + \La_i, \beta) $ for $\beta \in \rlQ$.
Note that there is sign difference in the notations between this paper and the previous papers \cite{KKOP21, KKOP22} (see Remark \ref{rem:changes}). 
Under the categorification, the determinantial module $\dM(w\La_i,\La_i)$ correspond to the \emph{unipotent quantum minor} $D(w\La_i, \La_i)$, which is a frozen variable of the quantum cluster algebra $A_q(\n(w))$. Thus the localization $\tcatC_w$ of $\catC_w$ with respect to the above left braiders categorifies the localization $A_q(\n(w))[  D(w\La_i, \La_i)^{-1} ; i\in I ]$ of $A_q(\n(w))$ at the frozen variables $D(w\La_i, \La_i)$ ($i\in I$). Note that, via the categorification, the Grothendieck ring $ K(\tcatC_{w})$ can be understood as a $q$-deformation of the coordinate ring of the unipotent cell $N^w$ associated with $w$ (see \cite[Theorem 4.13]{KO21}).
It turns out that the natural inclusion functor $ \iota_w\cl\catC_w \rightarrowtail R\gmod$ induces an equivalence of categories
$$
\widetilde \iota_w \cl  \tcatC_w \buildrel \sim \over \longrightarrow \bl R\gmod \br[\dM(w\La_i,\La_i)^{\circ -1}; i\in I].
$$ 
The localized category $\tcatC_w$ is rigid (\cite{KKOP21} for the left rigidity and \cite{KKOP22} for the full rigidity) and the dual functor on $\tcatC_w$ corresponds to the quantum twist automorphism (\cite{KO21}) via the categorification.

The subcategory $\catC_{w,v}$ was introduced in \cite{KKOP18} for a pair of Weyl group elements $w$, $v$ with $ w \ge v$ in the Bruhat order.
The category $\catC_{w,v}$ is defined as the intersection of two categories $\catC_w$ and $\catC_{*, v}$ (see Section \ref{Sec: Cwv} for definitions), and when $v = \id $, the category $\catC_{w,v}$ coincides with $\catC_w$.  The category $\catC_{w,v}$ is monoidal and the Grothendieck ring $K(\catC_{w,v})$ can be viewed as a $q$-deformation of the doubly-invariant algebra $^{N'(w)} \C[N] ^{N(v)}  $ (see \cite[Remark 2.19]{KKOP18} for precise notations).  The localization of $^{N'(w)} \C[N] ^{N(v)}  $ at the unipotent minors $D(w\La_i, v\La_i)$ is isomorphic to the coordinate ring $\C[R_{w,v}]$ of the \emph{open Richardson variety} $R_{w,v}$ and  it has a cluster algebra structure (see \cite{Lec16, CGGLSS22, GLSB23} and also \cite{CK22, LY19, Mer22, SSB22}). 
Under the categorification, the unipotent minors $D(w\La_i, v\La_i)$ correspond to the determinantial modules $\dM(w\La_i, v\La_i)$ which are central (up to a grading shift) in the category $\catC_{w,v}$.  Like the case of $\catC_w$, it is proved in \cite[Section 4]{KKOP22} that  there exist the graded left braiders in the category $\catC_{*, v}$
$$
( \dM(w\La_i, v\La_i),  \coR_{\dM(w\La_i, v\La_i)}, \dphi^l_{ w, v, \La_i} ) \qquad \text{ for any $i\in I$ }
$$ 
using the homomorphism induced from R-matrices
$$
\coR_{\dM(w\La_i, v\La_i)}(X) \cl \dM(w\La_i, v\La_i) \conv  X  \longrightarrow q^{- \dphi^l_{ w,v, \La_i} ( \wt(X) ) }  X  \conv   \dM(w\La_i,v\La_i),  
$$
for any $X \in \catC_{*, v}$, and  $\dphi^l_{w,v, \La_i} ( \beta) = - (w \La_i + v\La_i, \beta) $ for $\beta \in \rlQ$. 
We thus obtain the localization $\tcatC_{w,v}$ of $\catC_{w,v}$ with respect to the above left braiders, and
the Grothendieck ring $K(\tcatC_{w,v}) $ can be understood as a $q$-deformation of the coordinate ring $\C[R_{w,v}]$ of the open Richardson variety $R_{w,v}$.
The natural inclusion functor $ \iota_{w,v}^l\cl\catC_{w,v} \rightarrowtail \catC_{*, v}$ induces an equivalence of categories
$$
{\widetilde \iota}_{w,v}^{\ms{5mu}l} \cl  \tcatC_{w,v} \isoto\; \bl \catC_{*, v} \br[\dM(w\La_i, v\La_i)^{\circ -1}; i\in I]
$$ 
(see \cite[Theorem 4.5]{KKOP22}).
It is conjectured that the localized category $\tcatC_{w,v}$ gives a monoidal categorification for the $q$-deformation of the coordinate ring $\C[R_{w,v}]$ as a quantum cluster algebra.

\smallskip

In this paper, we develop a localization procedure by \emph{right braiders} and construct the localization $\tcatC_{w,v}$ of the category $\catC_{w,v}$ using the right braiders arising from  determinantial modules.
We then prove interesting categorical properties of the localized category $\tcatC_{w,v} $ including the right rigidity.

Let $\catT =\soplus_{\la \in\La} \catT_\la$ be a graded monoidal category with good properties (see Section \ref{Sec: right braider} for details). 
A triple $(C, \coRr_C, \phi)$ is called a \emph{graded right braider} in $\catT$ if $C$ is an object, $\phi\cl \La \to \Z$ is a $\Z$-linear map, and  
\begin{align} \label{Eq: rb}
\coRr_{C} (X) \cl X \tens C \to q^{-\phi(\la)}\tens C\tens X
\end{align}
is a morphism functorial in $X \in \catT_\la$ such that they satisfy a certain compatibility (see \eqref{eq:right braiders}). Following the same arguments in the theory of localization using left braiders,  we define the notion of a \emph{real commuting family of graded right braiders}, and show that the category $\catT$ can be localized by right braiders. 
The localization by the right braiders has the same properties as those of the localization by left braiders (see Theorem \ref{Thm: graded localization}).  
The difference between right and left braiders is the order of tensor products  in the morphisms $\coRr_{C}$ and $\coRl_{C}$ (see \eqref{Eq: rb} and \eqref{Eq: lb}).  
Although the localization procedure by right braiders  is a novel and crucial step in the paper, we skip several details in the section for localization by right braiders because the arguments used in the right braider setting are 
same as those used in the left braider setting (see Section \ref{Sec: right braider}).

The localization by right braiders is applied to the category $\catC_{w,v}$. We first prove that, for any dominant weight $\La$, there exists a morphism in $\catC_w$ 
 $$
 \coRr_{\dM(w\La, v\La)}(X) \cl X \conv \dM(w\La,v\La) \to q^{-(\wt(X),w\La+v\La)} \dM(w\La,v\La ) \conv X
 $$ 
functorial  in $X\in \catC_w$ (see Proposition \ref{prop:right_braider}). Since the determinantial module  $\dM(w\La,v\La)$ with $ \coRr_{\dM(w\La, v\La)}$ gives us a real commuting family of graded right braiders 
$$
\{\bl \dM(w\La_i,v\La_i), \coRr_{\dM(w\La_i,v\La_i)},\phi^r_{w,v,\La_i}\br \}_{i\in I}
$$ 
in the category  $\catC_w$ (see Corollary \ref{cor:right_baraider}), we can localize the categories $\catC_w$ and $\catC_{w,v}$ by the above right braiders, which are denoted by $\catC_w[\dM(w\La_i,v\La_i)^{\circ -1}; i\in I]$ and $ \catC_{w,v}[\dM(w\La_i,v\La_i)^{\circ -1}; i\in I]$. Since $\dM(w\La,v\La)$ is central in the category $\catC_{w,v}$ (up to a grading shift), the localization of $\catC_{w,v}$ by the right braiders coincides with that by the left braiders, i.e.,  
$$
\tcatC_{w,v} = \catC_{w,v}[\dM(w\La_i,v\La_i)^{\circ -1}; i\in I].
$$
The natural inclusion functor $ \iota_{w,v}\cl \catC_{w,v} \rightarrowtail \catC_{w}$ induces an equivalence of categories
$$
\widetilde \iota_{w,v} \cl  \tcatC_{w,v} \buildrel \sim \over \longrightarrow  \catC_{w} [\dM(w\La_i, v\La_i)^{\circ -1}; i\in I]
$$ 
(see Theorem \ref{thm:essential_surj}), and the kernels of 
the localization functors
$$
\Ql_{w,v}\cl\catC_{*,v}\to\tcatC_{w,v} \quad  \text{and} \quad  \Qr_{w,v}\cl\catC_{w}\to\tcatC_{w,v}
$$
with respect to the left and right  braiders respectively are characterized in terms of the $\Z$-valued invariant $\La$ (see Proposition \ref{prop:Laneq1} and Proposition \ref{prop:Laneq2}).
In the course of proofs, the membership criteria for $\catC_{*, v}$ and $\catC_{w,v}$ (Theorem \ref{thm:C*v}, Corollary \ref{coro:belongtoCwv1} and Corollary \ref{coro:belongtoCwv2}) are used in a crucial way.

We step further into the category $\catC_{w,v}$. We  prove that the localized category $\tcatC_{w,v}$ is right rigid by using the rigidity of 
$\tcatC_w$ 
(see Theorem \ref{thm:right_rigidity}).
It is interesting to ask how the right dual functor of $\tcatC_{w,v}$ is related to the twist automorphism for the Richardson variety $R_{w,v}$ introduced in \cite{GL22, LY19}.  
 It is  conjectured that $\tcatC_{w,v}$ is a \emph{rigid} monoidal category.     
We next interpret the $T$-systems consisting of quantum unipotent minors $D(w\La_i, v\La_i)$ at the categorical level to obtain the corresponding short exact sequences in terms of determinantial modules in $\catC_{w,v}$ (see Proposition \ref{prop:lamu=Lai}). Theorem \ref{thm: equiv between tcwv} says that there is an equivalence of monoidal categories 
$$
\tcatC_{w,v} \simeq \tcatC_{ws_i,vs_i}
$$
for a pair of Weyl group elements $w$ and $v$ such that  $w \ge v$, $ws_i > w$  and $vs_i >v$.

\smallskip
This paper is organized as follows.
In Section \ref{Sec: Preliminaries}, we review some preliminaries on quiver Hecke algebras and explain briefly the localization procedure by right braiders.
In Section \ref{Sec: loc rb}, the localization procedure is applied to the category  $\catC_{w,v}$ by studying the right braiders arising from determinantial modules. 
In Section \ref{Sec: properties of Cwv}, we prove interesting categorical properties of the localized category $\tcatC_{w,v} $ including the right rigidity of $\tcatC_{w,v}$.

\vskip 2em

\section{Preliminaries} \label{Sec: Preliminaries}
\subsection{Localization of a monoidal category via right braiders} \label{Sec: right braider}
In this subsection we recall the localization of a monoidal category via a real commuting family of  \emph{right  braiders} following \cite[Section 2]{KKOP21}. Since the localization via the \emph{left braiders} is studied in detail 
 in \cite[Section 2]{KKOP21},
we recall the right braiders case without proofs.

Let $\cor$ be a commutative ring,  $\La$ an abelian monoid,  and $(\catT, \tens)$  a $\cor$-linear monoidal category with a unit object $\one$ (see \cite[Section 1.5]{KKOP21} for the definition of monoidal categories and related notions).
 Assume that there is a direct sum decomposition of the category $\catT=\soplus_{\la \in\La} \catT_\la$ such that  $\tens$ induces a bifunctor $\catT_\la \times \catT_\mu \to \catT_{\la+\mu}$ for $\la,\mu\in \La$ and $\one \in \catT_0$.  We call $\catT$ a \emph{$\La$-graded monoidal category}.
 Let $q$ be an invertible central object in  $\catT_0$.   We write $q^n$ ($n\in\Z$) for $q^{\tens n}$ for the sake of simplicity. 
We say that  $f\in \Hom_{\catT}(q^d X,Y )$ is a morphism from $X$ to $Y$ of degree $d$,  and write  $d=\deg(f)$.

\begin{definition}
A \emph{graded right braider in $\catT$} is a triple $(C, \coRr_C, \phi)$ of an object $C$, a $\Z$-linear map $\phi\cl \La \to \Z$, and a morphism functorial in $X \in \catT_\la$
\eqn
\coRr_{C} (X) \cl X \tens C \to q^{-\phi(\la)}\tens C\tens X
\eneqn 
such that the following diagrams commute for any $X \in \catT_\la$ and $Y \in \catT_\mu$: 
\begin{equation} \label{eq:right braiders}
\begin{aligned}
\xymatrix{
X \otimes Y \otimes C  \ar[rr]^{X\otimes \coRr_C(Y)}  \ar[drr]_{\coRr_C(X \otimes Y)\ \ }  &  &q^{-\phi(\mu)} \tens X \otimes C \otimes Y  \ar[d]^{ \coRr_C(X)\otimes Y}   \\
& & q^{-\phi(\mu+\la)} \tens C\otimes X \otimes Y,
}
\ \qquad
\xymatrix{
\triv \otimes C   \ar[rr]^{\coRr_C(\triv)}  \ar[drr]_{ \simeq }  &  &   C  \otimes \triv \  \ar[d]^{ \bwr}   \\
& &    C.
}
\end{aligned}
\end{equation}
\end{definition}

A graded right braider $(C, \coRr_{C}, \phi)$ is called a \emph{central object} if $\coRr_{C}(X)$ is an isomorphism for any $X\in \catT$.
 
Let us denote the category of graded right braiders by $\catT^r_{br}$.  Note that $\catT^r_{br}$ is a monoidal category and there is a canonical faithful monoidal functor $\catT^r_{br} \to \catT$.

\begin{definition}
Let $I$ be an index set.  A family of graded right braiders $\{(C_i, \coRr_{C_i},\phi_i)\}_{i\in I}$ is called a \emph{real commuting family of graded right braiders in $\catT$} if
\bna
\item  $C_i\in \catT_{\la_i}$ for some $\la_i\in \La$,  and  $\phi_i(\la_j)+\phi_j(\la_i)=0$ for any $i,j \in I$, 
\item  $\coRr_{C_i}(C_i)  \in \cor^\times \id_{C_i\tens C_i}$ for any $i \in I$,
\item $\coRr_{C_j}(C_i) \circ \coRr_{C_i}(C_j)  \in \cor^\times \id_{C_j\tens C_i}$ for any $i,j \in I$.
\ee
\end{definition}

Define a $\Z$-linear map 
\eqn 
\phi \cl  \Z^{\oplus I} \times \La \to \Z \quad \text{given by} \quad (e_i, \la) \mapsto \phi_i(\la),
\eneqn
where $\{e_i\}_{i\in I}$ denotes the standard basis of $\Z^{\oplus I}$. 
We denote by $\phi_\al$ the $\Z$-linear map
 \eqn 
 \phi_\al\seteq\phi(\al,-)\cl  \La \to \Z \quad \text{ for each} \ \al \in \Z^{\oplus I}.
 \eneqn
 Note that one can choose a $\Z$-bilinear map $H \cl  \Z^{\oplus I} \times \Z^{\oplus I} \to \Z $ such that
$$\phi_i(\la_j) = H(e_i,e_j)-H(e_j,e_i) \quad \text{ for any} \ i,j\in I.$$ 
Then we have
\eqn
\phi(\al, L(\beta)) = H(\al,\beta)-H(\beta,\alpha) \quad \text{for any } \ \al,\beta\in \Z^{\oplus I},
\eneqn
where $L\cl \Z^{\oplus I}  \to \La$ be the $\Z$-linear map given by $e_i \mapsto \la_i$ for $i\in I$.

\begin{lem}[{\cite[Lemma 2.3, Lemma 1.16]{KKOP21}}]
Let $\st{( C_i ,  \coRr_{C_i}, \phi_i )}_{i\in I}$
be a real commuting family of right graded braiders in  $\catT$.
\bnum
\item
There exists a family $\{\eta_{ij}\}_{i,j\in I}$
of elements in $\corp^\times$ such that
\eqn
\coRr_{C_i}(C_i)&&=\eta_{ii}\; \id_{C_i \otimes C_i},\\
\coRr_{C_j}(C_i) \circ \coRr_{C_i}(C_j)&& = \eta_{ij}\eta_{ji}\;\id_{C_j \otimes C_i}
\eneqn
for all $ i ,j \in I$.

\item 
There exist
a graded right braider $C^\al=(C^\alpha, \coRr_{C^\alpha}, \phi_\alpha )$ for each $\alpha \in \Z^{\oplus I}_{\ge0}$,
and an isomorphism $\xi_{\alpha, \beta}\cl  C^\alpha \otimes C^\beta \buildrel \sim\over \longrightarrow q^{H(\al,\beta)}\tens C^{\alpha+\beta}$ in $\catT^r_{br}$ for $\alpha, \beta \in \Z^{\oplus I}_{\ge0}$
such that 
\bna
\item  $C^0=\triv$  and $ C^{e_i} = C_i $ for $i\in I$,
\item  the diagram in $\catT^r_{br}$ 
\begin{equation} \label{Eq: xi sum}
\begin{aligned}
\xymatrix{
C^\alpha \otimes C^\beta \otimes C^\gamma \ar[d]_{ C^\alpha \otimes \xi_{\beta, \gamma}} \ar[rr]^{\xi_{\alpha, \beta} \otimes C^\gamma} && \ar[d]^{\xi_{\alpha+\beta, \gamma}} q^{H(\al,\beta)}\tens C^{\alpha+\beta}\otimes C^\gamma \\
q^{H(\beta,\gamma)} \tens C^\alpha \otimes C^{\beta + \gamma} \ar[rr]^{\xi_{\alpha, \beta+\gamma}} && q^{H(\al,\beta)+H(\al,\gamma)+H(\beta,\gamma)} \tens C^{\alpha+\beta +\gamma}
}
\end{aligned}
\end{equation}
commutes for any $\alpha, \beta, \gamma \in \Z^{\oplus I}_{\ge0}$,
\item the diagrams  in $\catT^r_{br}$ 
\begin{equation} \label{Eq: CF ij}
\begin{aligned}
\ba{ccc}
\xymatrix{
C^0 \otimes C^0  \ar[d]^-{ \bwr } \ar[rr]^{\xi_{0,0} } && \ar[d]^-{ \bwr } C^{0}\\
\triv \otimes \triv  \ar[rr]^{ \simeq} && \ \ \triv  \ ,
}
\ba{c}\\[3ex]\qtq\ea\quad
\xymatrix{
C^\alpha \otimes C^\beta \ar[d]_{ \xi_{\alpha, \beta}} \ar[rr]^{\coRr_{C^\beta}(C^\alpha) } && \ar[d]^{\xi_{ \beta, \alpha}}q^{-\phi(\beta,L(\al))}\tens  C^\beta \otimes C^\alpha \\
q^{\gH(\al,\beta)}\tens C^{\alpha + \beta} \ar[rr]^{\bce(\alpha,\beta)\; \id_{C^{\alpha+\beta}} } && q^{H(\al,\beta)}\tens C^{\alpha+\beta}
}
\ea
\end{aligned}
\end{equation}

commute for any $i,j\in I$ and $\alpha, \beta, \gamma \in \Z^{\oplus I}_{\ge0}$, where
\begin{align} \label{Eq: eta}
\bce(\alpha, \beta) \seteq  \prod_{i,j \in I} \bce_{i,j}^{a_ib_j} \in \bR^\times
\quad \text{for $\alpha = \sum_{i\in I} a_i e_i$ and $\beta = \sum_{j\in I} b_j e_j$ in $\Z^{\oplus I}$.}
\end{align}

\end{enumerate}
\end{enumerate}

\end{lem}
Define a partial order $\preceq$ on $ \Z^{\oplus I}$ by
$$
\alpha \preceq \beta  \quad \text{ for } \alpha, \beta \in \Z^{\oplus I} \text{ with }  \beta - \alpha \in \Z^{\oplus I}_{\ge0},
$$
and set
$$
\Ds_{\alpha, \beta} \seteq  \{ \delta \in \Z^{\oplus I} \mid \alpha + \delta , \beta+\delta \in \Z^{\oplus I}\}
$$
for $\alpha,\beta \in \Z^{\oplus I}$.

For $X \in \catT_\lambda$,  $Y \in \catT_\mu$  and $ \delta  \in \Ds_{\alpha, \beta}$,
we set
\eq
\gHm_\delta( (X, \alpha  ), (Y, \beta  ) ) \seteq  \Hom_{\catT}(  X\otimes C^{\delta + \alpha}, q^{\gH( \beta-\alpha,\delta)-\dphi(\delta+\beta, \mu)} \otimes   C^{ \delta + \beta} \otimes Y).
\eneq
For $ \delta, \delta' \in \Ds_{\alpha, \beta}$ with $\delta \preceq  \delta'$ and $ f \in \gHm_\delta( (X, \alpha  ), (Y, \beta  ) )$,
we define $\gzeta_{\delta', \delta}(f) \in \gHm_{\delta'}( (X, \alpha  ), (Y, \beta  ) )$  to be the morphism such that the following diagram commutes:
$$
\xymatrix@C=4em{
X  \tens C^{ \delta + \alpha} \tens C^{\delta' - \delta}     \ar[dd]^{\bwr}_{  \xi_{ \delta+\alpha,\delta'-\delta} }  \ar[rr]^{ f \tens C^{\delta'-\delta}  \qquad \qquad }   &&
q^{\gH(\beta-\alpha,\delta) -\dphi(\delta+\beta, \mu)} \otimes    C^{ \delta+\beta}  \otimes Y \otimes C^{\delta' - \delta}   \ar[d]^{  \coRr_{C^{\delta'-\delta}} (Y) } \\
 &&  q^{\gH(\beta-\alpha,\delta) -\dphi(\delta'+\beta, \mu)}  \otimes C^{ \delta + \beta} \otimes C^{\delta' - \delta} \otimes Y
  \ar[d]_{\bwr}^{\xi_{ \delta+\beta , \delta'-\delta}}  \\
q^{\gH( \delta+\alpha,\delta'-\delta)} \otimes X \otimes C^{  \delta'+\alpha}   \ar[rr]^{ q^{\gH( \delta+\alpha,\delta'-\delta)} \tens \gzeta_{\delta', \delta}(f) \qquad \qquad} &&
q^{\gH(\beta-\alpha,\delta) -\dphi(\delta'+\beta, \mu) + \gH(\delta+\beta, \delta'-\delta) } \otimes C^{ \delta'+\beta}  \otimes Y.
}
$$
Then,
$\gzeta_{\delta', \delta } $ is a map from $  \gHm_\delta( (X, \alpha  ), (Y, \beta  ) )$  to $ \gHm_{\delta'}( (X, \alpha  ), (Y, \beta  ) )$ and   $\zeta^\gr_{\delta'', \delta' }  \circ \zeta^\gr_{\delta', \delta } = \zeta^\gr_{\delta'', \delta } $
for $\delta \preceq \delta' \preceq \delta''$, so that $\{\zeta^\gr_{\delta', \delta} \}_{\delta,\delta' \in \Ds_{\al,\beta}}$ forms an inductive system indexed by the set $\Ds_{\al,\beta}$.

Define 
\begin{align*}
\Ob (\lT) &\seteq  \Ob(\catT) \times \Z^{\oplus I}
\end{align*}
and for  $X \in \catT_{\lambda}$ and $Y \in \catT_{\mu}$ define
\begin{align*}
\Hom_{\lT}( (X, \alpha), (Y, \beta) ) &\seteq   \indlim[{  \substack{\delta \in \Ds_{\alpha, \beta}, \\  \lambda + \gL(\alpha) = \mu + \gL(\beta)  }    }]
\Hm^\gr_{\delta}( (X, \alpha  ), (Y, \beta  ) ).
\end{align*}

Let $X \in \catT_\lambda$, $Y \in \catT_\mu$ and $Z \in \catT_\nu$. For  $f \in \gHm_\delta( (X, \alpha  ), (Y, \beta  ) ) $ and  $ g \in  \gHm_\epsilon( (Y, \beta  ), (Z, \gamma  ) )$, we define
$$
\gPsi_{\delta, \epsilon} (f,g)  \seteq  \bce(\delta+ \beta, \beta-\gamma) \cdot \gtPsi_{\delta, \epsilon}(f,g) ,
$$
where  $\gtPsi_{\delta, \epsilon}(f,g)$ is the morphism such that the following diagram commutes:
$$
\xymatrix{
X\otimes C^{\delta+\alpha} \otimes C^{\epsilon+\beta} \ar[dd]^{\wr}_{\xi_{\delta+\alpha, \epsilon+\beta} }  \ar[rr]^f &&  q^{a}\otimes  C^{\delta+\beta}\otimes   Y \otimes C^{\epsilon+\beta} \ar[d]^g \\
&& q^{b}\otimes  C^{\delta+\beta} \otimes C^{\epsilon+\gamma}\otimes Z   \ar[d]_{\wr}^{\xi_{ \delta+\beta,\epsilon+\gamma}} \\
q^{\gH(\delta+\alpha,\epsilon+\beta)}\otimes   X\otimes  C^{\delta+\epsilon+\beta+\alpha} \ar[rr]^{\qquad \gtPsi_{\delta, \epsilon}(f,g)} && q^{c}\otimes C^{\delta+\epsilon+\beta+\gamma} \otimes  Z  ,
}
$$
where
\begin{align*}
&a = \gH(\beta-\alpha.\delta)-\dphi(\delta+\beta, \mu), \quad b = a+\gH(\gamma-\beta,\epsilon) -\dphi(\epsilon+\gamma, \nu), \\
& c = b + \gH(\delta+\beta,\epsilon+\gamma).
\end{align*}
We may assume that $\mu=L(\gamma-\beta)+\nu$.  Hence we  have
\begin{align*}
c -\gH(\delta+\alpha,\epsilon+\beta)= \gH ( \gamma-\alpha,\beta+\epsilon+\delta )-\dphi(\delta+\epsilon+\beta+\gamma,\nu)
\end{align*}
so that
$$
\gPsi_{\delta, \epsilon} (f,g) \in \gHm_{\delta+\epsilon+\beta}( (X, \alpha  ), (Z, \gamma ) ) .
$$
It follows that 
\begin{align*}
\gPsi_{ \delta', \epsilon'}( \gzeta_{\delta', \delta}(f), \gzeta_{\epsilon', \epsilon}(g) ) = \gzeta_{  \delta' + \epsilon' + \beta, \delta+\epsilon+\beta }(\Psi_{\delta, \epsilon}(f,g)),
\end{align*}
which yields the composition in $\lT$:
$$
\Hom_{\lT}( (X, \alpha), (Y, \beta) )  \times \Hom_{\lT}( (Y, \beta), (Z, \gamma) )  \rightarrow
\Hom_{\lT}( (X, \alpha), (Z, \gamma) ).
$$

Because this  composition in $\lT$ is associative,  $\lT$ becomes a category.
 By the construction, we have the decomposition
$$
\lT = \bigoplus_{\mu \in \Lambda} \lT_{\mu}, \qquad \text{where }  \lT_{\mu} \seteq  \{ (X, \alpha) \mid X \in \catT_\lambda, \ \lambda + \gL(\alpha)=\mu  \}.
$$

The category $\lT$  is a monoidal category with the following  tensor product. For $\alpha, \alpha', \beta, \beta' \in \lG $, $X \in \catT_\lambda$, $X' \in \catT_{\lambda'}$, $Y \in \catT_\mu$ and $Y' \in \catT_{\mu'}$,
we define
$$
(X, \alpha) \otimes (Y, \beta) \seteq  ( q^{ \dphi(\alpha, \mu) + \gH(\alpha, \beta)} \otimes X \otimes Y, \alpha+\beta  ),
$$
and, for $f \in \gHm_\delta((X, \alpha), (X', \alpha'))$ and $g \in \gHm_\epsilon((Y, \beta), (Y', \beta'))$, we define
$$
\gT_{\delta, \epsilon}(f,g) \seteq  \eta(\delta, \beta-\beta') \gtT_{\delta, \epsilon}(f,g) ,
$$
where $ \gtT_{\delta, \epsilon}(f,g)$ is the morphism such that the following diagram commutes:
$$
\xymatrix{
X\otimes C^{\delta+\alpha}  \otimes Y \otimes C^{\epsilon+\beta} \ar[rr]^{f \otimes g} &&  q^{b}\otimes  C^{\delta+\alpha'} \otimes X' \otimes  C^{\epsilon+\beta'}  \otimes  Y'    \ar[d]^{\coRr_{C^{\epsilon+\beta'}} (X') } \\
q^{\dphi(\delta+\alpha, \mu) }\otimes  X \otimes Y  \otimes C^{\delta+\alpha} \otimes C^{\epsilon+\beta} \ar[u]^{\coRr_{C^{\delta+\alpha}} (Y) }  \ar[d]^\bwr_{\xi_{\delta+\alpha, \epsilon+\beta}}
&&  q^{c}\otimes  C^{\delta+\alpha'} \otimes C^{\epsilon+\beta'}  \otimes X' \otimes Y'  \ar[d]_\bwr^{\xi_{\delta+\alpha', \epsilon+\beta'}} \\
q^{a}\otimes  X \otimes Y \otimes C^{\delta+\epsilon+\alpha+\beta} \ar[rr]^{\gtT_{\delta,\epsilon} (f, g)}  &&  q^{d} \otimes C^{\delta+\epsilon+\alpha'+\beta'} \otimes  X' \otimes Y',
}
$$
for some $a,b,c,d\in\Z$. 

Note that we may assume that $\mu+L(\beta)=\mu'+L(\beta')$.
Hence we have
\begin{align*}
&d-a = \dphi(\al',\mu')+\gH(\alpha',\beta') -\left( \dphi(\alpha,\mu)+ \gH(\alpha,\beta)\right) \\
&\qquad + \gH(\alpha'+\beta'-\alpha-\beta, \delta+\epsilon) - \dphi(\delta+\epsilon+\alpha'+\beta',\la'+\mu').
\end{align*}
It follows that
$$
\gT_{\delta, \epsilon}(f,g) \in \gHm_{\delta+\epsilon}((X, \alpha)\tens (Y,\beta), (X', \alpha') \otimes (Y', \beta') ).
$$
Then we have
\eqn
\gT_{\delta', \epsilon'}(\gzeta_{\delta',\delta}(f), \gzeta_{\epsilon',\epsilon}(g)) =
\gzeta_{\delta'+\epsilon',\delta+\epsilon} (\gT_{\delta,\epsilon}(f,g))
\quad \text{for $\delta' \succeq \delta$ and $\epsilon' \succeq \epsilon$.}
\eneqn
That is,  the map $\gT_{\delta, \epsilon}$ is compatible with the inductive system and hence it induces a well-defined map
\begin{align} \label{Eq: tensor for localization}
f \otimes g \in \Hom_{\lT}( (X, \alpha) \otimes (Y, \beta) , (X', \alpha') \otimes (Y', \beta') )
\end{align}
for
$f \in  \Hom_{\lT}( (X, \alpha)  , (X', \alpha')  )$ and $g \in  \Hom_{\lT}( (Y, \beta)  , (Y', \beta')  )$.

Moreover,  we have
$$
\Psi_{\delta_1+\delta_2, \epsilon_1 + \epsilon_2} (T_{\delta_1, \delta_2} (f_1, f_2),  T_{\epsilon_1, \epsilon_2} (g_1, g_2) )
= T_{\delta_1+\epsilon_1+\beta_1, \delta_2+\epsilon_2 + \beta_2 } (\Psi_{\delta_1, \epsilon_1} (f_1,g_1), \Psi_{\delta_2, \epsilon_2} (f_2,g_2) )
$$
where $f_k \in \Hm_{\delta_k}( (X_k, \alpha_k), (Y_k, \beta_k) )  $ and $g_k \in \Hm_{\epsilon_k}( (Y_k, \beta_k), (Z_k, \gamma_k) ) $ for $k=1,2$ (see \cite[Proposition 2.5]{KKOP21}).

It follows that the map $\tens$ on $\lT$ defines a bifunctor $\tens \cl  \lT \times \lT \to \lT$.

\begin{thm}  \label{Thm: graded localization}
Let $ \st{C_i=(C_i, \coRr_{C_i}, \dphi_i )}_{i\in I} $ be a 
real commuting family of graded right braiders  in $\catT$. 
Then the category
$\lT$ defined above becomes a monoidal category. 
There exists a monoidal functor 
$\Upsilon\cl \catT \to \lT$ 
and a real commuting family of graded  right braiders $\st{\tC_i=(\tC_i, \coRr_{\tC_i}, \phi_i)}_{i\in I}$ in $\lT$ satisfy the following properties:

\bnum
\item
for $i\in I$, $\Upsilon(C_i) $ is isomorphic to $ \widetilde{C}_i$ and it is invertible in $(\lT)_{\,\mathrm{br}}$, 
\item 
for $i\in I$ and $X\in\catT_\la$, the diagram
$$
\xymatrix{
\Upsilon(X \otimes C_i )  \ar[r]^\sim  \ar[d]_{\Upsilon( \coRr_{C_i} (X)  )\ms{10mu}}^-\bwr  & \Upsilon(X) \otimes  \widetilde{C}_i  \ar[d]_{ \coRr_{\widetilde{C}_i} (\Upsilon(X)  )\ms{5mu}} ^-\bwr \\
\Upsilon(   q^{ - \phi_i(\la) }   \otimes C_i \otimes X )  \ar[r]^\sim 
& q^{ - \phi_i(\la) }\tens  \ \widetilde{C}_i \otimes \Upsilon(X)
}
$$
commutes.
\setcounter{myc}{\value{enumi}}

\end{enumerate}

Moreover, the functor $\Upsilon$ satisfies the following universal property:
\bnum\setcounter{enumi}{\value{myc}}
\item  If there are another $\La$-graded monoidal category $\catT'$ 
with an invertible central object $q\in\catT'_0$
with and a $\La$-graded  monoidal functor $\Upsilon'\cl  \catT \rightarrow \catT'$ 
such that  
\bna
\item  $\Upsilon'$ sends the central object $q\in\catT_0$ to $q\in\catT'_0$, 
\item   \label{Eq: loc 1}
$\Upsilon'(C_i) $ is invertible in $\catT'$ for any $i\in I$ and
\item 
for any $i\in I$ and $X\in\catT$, $\Upsilon (\coRr_{C_i}(X))\cl
\Upsilon'(X\tens C_i)\to\Upsilon'(q^{  - \phi_i(\la) }\tens  C_i\tens X)$ is an isomorphism,
\end{enumerate}
then there exists a monoidal functor $\mathcal F$, which is unique up to a unique isomorphism,  such that
the diagram
$$
\xymatrix{
\catT \ar[r]^{\Upsilon} \ar[dr]_{\Upsilon'}  & \lT \ar@{.>}[d]^{\mathcal F }\\
& \catT'
}
$$
commutes.
\end{enumerate}

\end{thm}

We denote by $\catT [ C_i^{\otimes -1} \mid i\in I]$ the localization $\lT$ in Theorem \ref{Thm: graded localization}.
If $\catT$ is an abelian monoidal category with exact tensor product,  then so is $\catT [ C_i^{\otimes -1} \mid i\in I]$,  and the functor  $\Upsilon\cl  \catT \rightarrow \catT [ C_i^{\otimes -1} \mid i\in I]$  is an exact monoidal functor. 

Note that
$$
(X, \alpha+\beta) \simeq (q^{-\gH(\alpha,\beta)}  X\tens C^\alpha, \beta), \quad
(\triv, \beta) \otimes  (\triv, -\beta) \simeq q^{-\gH(\beta, \beta)} (\triv, 0)
$$ 
for $\alpha \in \Z^{\oplus I}_{\ge 0}$ and  $\beta \in  \Z^{\oplus I}$.

\begin{remark} \label{rem:changes}
\begin{enumerate}
\item
Recall that a \emph{graded left braider in $\catT$} is a triple $(C, \coR_C, \phi)$ of an object $C$, a $\Z$-linear map $\phi\cl  \La \to \Z$, and a morphism functorial in $X \in \catT_\la$
\eqn
\coR_{C} (X) \cl  C \tens X  \to q^{-\phi(\la)} X\tens C
\eneqn 
with analogous conditions to  \eqref{eq:right braiders}. 
The materials in this subsection, containing the above theorem, are proved in \cite[Section 2]{KKOP21} for the localization via a real commuting family of graded \emph{left braiders}.  
\item  We changed $\phi$  in \cite[Section 2.3]{KKOP21} with $-\phi$ in this paper in order to have $\phi(\la)=\deg(\coRr_C(X))$ when $X \in \catT_\la$. 
\end{enumerate}
\end{remark}

\subsection{Quiver Hecke algebras }
A Cartan datum $ \bl\cmA,\wlP,\Pi,\Pi^\vee,(\cdot,\cdot) \br $ is a quintuple of a generalized Cartan matrix $\cartan$,    a free abelian group $\wlP$,  a set of simple roots, $\Pi = \{ \alpha_i \mid i\in I \} \subset \wlP$,  a  
 set of simple coroots $\Pi^{\vee} = \{ h_i \mid i\in I \} \subset \wlP^{\vee}\seteq\Hom( \wlP, \Z )$ ,  and a $\Q$-valued 
symmetric bilinear form  $(\cdot,\cdot)$ on $\wlP$ such that
\bna 
\item $\cmA = (\langle h_i,\alpha_j\rangle)_{i,j\in I}$,
\item  $(\alpha_i,\alpha_i)\in 2\Z_{>0}$ for any $i\in I$,
\item $\langle h_i, \lambda \rangle =\dfrac{2(\alpha_i,\lambda)}{(\alpha_i,\alpha_i)}$ for $i\in I$ and $\lambda \in \Po$,
\item for each $i\in I$, there exists $\Lambda_i \in \wlP$
such that $\langle h_j, \Lambda_i \rangle = \delta_{ij}$ for any $j\in I$.
\end{enumerate}

Let $\rlQ\seteq\soplus_{i\in I} \Z\al_i$
and $\rlQ_+\seteq\soplus_{i\in I} \Z_{\ge0} \al_i$ be the root lattice and the positive root lattice of the symmetrizable Kac-Moody algebra $\g(\cartan)$,  respectively.

Let $\weyl$ be the Weyl group of $\g(\cartan)$,  the subgroup of
$\Aut(\wtl)$ generated by the simple reflections $\st{s_i}_{i\in I}$ where
$s_i(\la)=\la-\ang{h_i,\la}\al_i$ for $\la \in \wtl$.

\medskip
Let $(\qQ_{i,j}(u,v) \in \bR[u,v])_{i,j\in I}$ be a family of polynomials 
such that 
\begin{align}
\qQ_{i,j}(u,v) =\bc
                   \sum\limits
_{p(\alpha_i , \alpha_i) + q(\alpha_j , \alpha_j) = -2(\alpha_i , \alpha_j) } t_{i,j;p,q} u^pv^q &
\text{if $i \ne j$,}\\[3ex]
0 & \text{if $i=j$,}
\ec\label{eq:Q}
\end{align}
where $t_{i,j;-a_{ij},0} \in  \bR^{\times}$ and
$\qQ_{i,j}(u,v)= \qQ_{j,i}(v,u) \quad \text{for all} \ i,j\in I.$
We set
\begin{align*}
\bQ_{i,j}(u,v,w)\seteq\dfrac{ \qQ_{i,j}(u,v)- \qQ_{i,j}(w,v)}{u-w}\in \bR[u,v,w].
\end{align*}

For $\beta\in \rlQ_+$,  the set
$I^\beta\seteq  \Bigl\{\nu=(\nu_1, \ldots, \nu_n ) \in I^n \bigm| \sum_{k=1}^n\alpha_{\nu_k} = \beta \Bigr\}$ is stable under the the symmetric group $\mathfrak{S}_n = \langle s_k \mid k=1, \ldots, n-1 \rangle$ action  given by place permutations.

The height of $\beta=\sum_{i\in I} b_i \al_i\in\rtl$ is given by $\height{\beta}\seteq\sum_{i\in I} |b_i|$.
\begin{df}
Let $\beta\in\rlQ_+$ with $\height{\beta}=n$.
The {\em quiver Hecke algebra} $R(\beta)$ associated with the Cartan datum $ \bl\cmA,\Pi,\wlP,\Pi^\vee,(\cdot,\cdot) \br $ and the family of polynomials $(\qQ_{i,j}(u,v))_{i,j\in I}$
is the $\bR$-algebra generated by
$$
\{e(\nu) \mid \nu \in I^\beta \}, \; \{x_k \mid 1 \le k \le n \},
 \; \text{and} \; \{\tau_l \mid 1 \le l \le n-1 \}
$$
subject to the defining relations:
\eqn
&& e(\nu) e(\nu') = \delta_{\nu,\nu'} e(\nu),\ \sum_{\nu \in I^{\beta}} e(\nu)=1,\
x_k e(\nu) =  e(\nu) x_k, \  x_k x_l = x_l x_k,\\
&& \tau_l e(\nu) = e(s_l(\nu)) \tau_l,\  \tau_k \tau_l = \tau_l \tau_k \text{ if } |k - l| > 1, \\[1ex]
&&  \tau_k^2 = \sum_{\nu\in I^\beta}\qQ_{\nu_k, \nu_{k+1}}(x_k, x_{k+1})e(\nu), \\[5pt]
&& \tau_k x_l - x_{s_k(l)} \tau_k =
\bl\delta(l=k+1)-\delta(l=k)\br
\sum_{\nu\in I^\beta,\ \nu_k=\nu_{k+1}}e(\nu),\\
&&\tau_{k+1} \tau_{k} \tau_{k+1} - \tau_{k} \tau_{k+1} \tau_{k}
=\sum_{\nu\in I^\beta,\ \nu_k=\nu_{k+2}}
\bQ_{\,\nu_k,\nu_{k+1}}(x_k,x_{k+1},x_{k+2}) e(\nu).
\eneqn
\end{df}

The algebra $R(\beta)$ is equipped with the $\Z$-grading given by
\eq
&&\deg(e(\nu))=0, \quad \deg(x_k e(\nu))= ( \alpha_{\nu_k} ,\alpha_{\nu_k}), \quad  \deg(\tau_l e(\nu))= -(\alpha_{\nu_{l}} , \alpha_{\nu_{l+1}}).
\eneq

We denote by $\Modg(R(\beta))$ 
the category of graded modules 
over $R(\beta)$ with grade preserving homomorphisms.
We denote by $q$ the grade shift functor: $(qM)_k=M_{k-1}$,
We write
$\HOM(M,N)=\soplus\HOM(M,N)_k$ with
$\HOM(M,N)_k=\Hom(q^kM,N)$. 

The full subcategory of $\Modg(R(\beta))$ consisting of the objects which are finite-dimensional  over $\cor$ is denoted by $R(\beta)\gmod$.
For $M\in R\gmod$,  the space $M^\star \seteq  \HOM_{\bR}(M, \bR)$ is  an  $R(\beta)$-module via  the graded $\cor$-algebra antiautomorphism of $R(\beta)$ which fixes the generators
 $e(\nu)$, $x_k$, and $\tau_k$'s. 
We say that $M$ is \emph{self-dual} if $M \simeq M^\star$ in $R\gmod$. 
 For each simple module $M$ in $R\gmod$,  there exists $m\in \Z$ such that $q^mM$
is self-dual.

For $\al,\beta\in\prtl$, we set
$$e(\al,\beta)=\sum_{\substack{\nu\in I^{\al+\beta}\\
\sum_{k=1}^{\height{\al}} \al_{ \nu_k} =\al,\ \sum_{k=1}^{\height{\beta}} \al_{ \nu_{k+\height{\al}}}
 =\beta}}\hs{-3ex}e(\nu)\hs{1ex}\in R(\al+\beta).$$

Then there is an injective $\cor$-algebra homomorphism $R(\al) \tens R(\gamma) \to e(\al,\beta) R(\al+\beta)e(\al,\beta)$ and hence we can define the \emph{convolution product} by
$$
M \conv N \seteq R(\al+\beta) e(\al, \beta) \otimes_{R(\al) \otimes R(\beta)} (M \otimes N)
$$
for
 $M \in \Modg(R(\al))$ and $N \in\Modg(R(\beta))$.

Set
$$\Modg(R) \seteq\soplus_{\la\in \rtl_-}\Modg(R)_{\la} \qt{and} \quad  \ R\gmod =\soplus_{\la\in \rtl_-}(R\gmod)_\la, $$
where
$$\Modg(R)_{\la}  := \Modg(R(-\la)) \qt{and}\quad
 R\gmod_{\la}  := R(-\la)\gmod \qt{for } \ \la \in \rtl_-.$$
Hence the categories $\Modg(R)$ and $R\gmod$ are  $\rtl_-$-graded monoidal categories.
We denote $\wt(X):=-\beta$ when $X \in \Modg(R)_{-\beta}=\Modg(R(\beta))$.

For $M,N\in R\gmod$, we denote by $M\hconv N$ the head of $M\conv N$ and
by $M\sconv N$ the socle of $M\conv N$.
\

 A simple module $M$ is  called \emph{real} if $M\conv M$ is simple.

\subsection{R-matrices} \label{subSec: R-matrices}
Let $\beta \in \rlQ_+$ with $m =  \Ht(\beta)$. For  $k=1, \ldots, m-1$ and $\nu \in I^\beta$, the \emph{intertwiner} $\varphi_k \in R(\beta) $ is defined by 
\eq
\varphi_k e(\nu) =
\bc
 \bl\tau_k(x_k-x_{k+1})+1\br e(\nu) 
& \text{ if } \nu_k = \nu_{k+1}, 
 \\
 \tau_k e(\nu) & \text{ otherwise.}
\ec\label{def:intertwiner}
\eneq
Since $\st{\vphi_k}_{1\le k<m}$ satisfies the braid relation, we can define
$\vphi_w$ for any element $w$ of the symmetric group $\sym_m$ of degree $m$.

For $m,n \in \Z_{\ge 0}$, we set $w[m,n]$ to be the element of $\sg_{m+n}$ such that
$$
w[m,n](k) \seteq
\left\{
\begin{array}{ll}
 k+n & \text{ if } 1 \le k \le m,  \\
 k-m & \text{ if } m < k \le m+n.
\end{array}
\right.
$$

Let $M\in \Modg\bl R(\beta)\br $ and $N\in \Modg\bl R(\gamma)\br$ and
define the $R(\beta)\otimes R(\gamma)$-linear map $M \otimes N \rightarrow N \conv M$  by $$u \otimes v \mapsto \varphi_{w[\height{\gamma},\height{\beta}]}(v \etens u).$$
Then it extends to an $R(\beta+\gamma)$-module homomorphism (neglecting a grading shift)
$$
\Runi_{M,N}\cl  M\conv N \longrightarrow N \conv M.
$$
For $\beta \in \rlQ_+$ and $i\in I$,  let $\mathfrak{p}_{i, \beta}$ be an element in the center $Z(R(\beta))$ of $R(\beta)$
\begin{align} \label{Eq: def of p}
\mathfrak{p}_{i, \beta}  \seteq \sum_{\nu \in I^\beta} \Bigl(\hs{1ex}  \prod_{a \in \{1, \ldots, \Ht(\beta) \},\ \nu_a=i} x_a \Bigr) e(\nu)\in Z(R(\beta)).
\end{align}

 Assume that $M$  is a simple module in $R(\beta)\gmod$,  and there exists  an $R(\beta)$-module $\Ma$ with an endomorphism $z_{\Ma}$ of $\Ma$
with degree $d_{\Ma} \in \Z_{>0}$ such that
\eq \label{eq:oldaff}
&&\hs{2ex}\parbox{75ex}{
\begin{enumerate}[\rm (i)]
\item $\Ma / z_{\Ma} \Ma \simeq M$,
\item $\Ma$ is a finitely generated free module over the polynomial ring $\bR[z_{\Ma}]$,
\item $\mathfrak{p}_{i,\beta} \Ma \ne 0$ for all $i\in I$.\label{it:nonzeroP}
\end{enumerate}
}\eneq
We call $(\Ma,z_{\Ma})$ an \emph{affinization of $M$}.

Let $\Ma$ be an affinization of a simple $R$-module $M$, and  let $N$ be a non-zero $R$-module. We define a homomorphism (up to a grading shift)
$$
\Rre_{\Ma, N} \seteq  \z^{-s} \Runi_{\Ma, N}\cl  \Ma \conv N \longrightarrow N \conv \Ma
$$
where $s$ is the largest integer such that $ \RR_{\Ma, N}(\Ma \conv N) \subset \z^s (N \conv \Ma)$.
Then the homomorphism (up to a grading shift)
$$
\rmat{M,N} \cl M \conv N \longrightarrow N \conv M
$$
 induced from $\nR_{\Ma, N}$ by specializing at $\z=0$ never vanishes.
 We call $\rmat{M,N}$ the \emph{r-matrix}  between $M$ and $N$.
Let
\begin{align*}
\La(M,N) \seteq  \deg (\rmat{M,N}) ,
\end{align*}
and define
\begin{align*}
\tLa(M,N) \seteq   \frac{1}{2} \bl \La(M,N) + (\wt(M), \wt(N)) \br ,\quad 
\Dd(M,N) \seteq  \frac{1}{2} \bl\La(M,N) + \La(N,M)\br.
\end{align*}
Note that  
$ \de(M,N)$ and $ \tLa(M,N) \in \Z_{\ge 0}$ are non-negative integers (\cite[Lemma 3.11]{KKOP21}).

\medskip
A real simple module which admits an affinization is called \emph{\afr}.
The following result is used frequently throughout the paper. 

\Prop [{\cite[Theorem 3.2]{KKKO15}, \cite[Proposition 3.2.9]{KKKO18}}]
 \label{prop:simplehd} 
Let $M$ and $N$ be simple modules in $R\gmod$. Assume that one of them is \afr.
Then, the convolution
$M\conv N$ has a simple head and a simple socle.
Moreover,  we have
$$\dim\HOM_{R\gmod}(M\conv N,N\conv M)=1$$ and
$$M\hconv N\simeq\Im(\rmat{M,N})\simeq N\sconv M
\qt{up to grading shifts.}$$
\enprop

\subsection{Partial order on the weight lattice}

We define the partial order $\ble$ on $\wtl$ as follows:
$\la \ble \mu$  for $\la,\mu\in\wtl$ if there exists a sequence  of positive real roots $\beta_1,\ldots,\beta_r$ such that $(\beta_k, s_{\beta_{k+1}} s_{\beta_{k+2}} \cdots s_{\beta_r}\mu) >0$ for all $1\le k\le r$ and $\la=s_{\beta_1} s_{\beta_2} \cdots s_{\beta_r} \mu$. 
We have $\mu-\la\in\prtl$ if $\la \ble \mu$. Hence $\ble$ is a partial order on $\wtl$.

\begin{lem} \label{lem:lamu} 
Let $\la \ble \mu$ and $\al$ be a simple root satisfying $(\al,\la) \le 0$.
Then we have 
\bnum
\item
if $(\al,\mu)\ge 0$, then we have $s_\al \la \ble \mu$,
\item if $(\al,\mu)\le0$, then we have $s_\al \la \ble s_\al \mu$.
\ee
\end{lem}
\begin{proof}
Let $\beta_1,\ldots,\beta_r$ be a sequence of positive real roots  such that $(\beta_k, s_{\beta_{k+1}} s_{\beta_{k+2}} \cdots s_{\beta_r}\mu) >0$ for all $1\le k\le r$ and $\la=s_{\beta_1} s_{\beta_2} \cdots s_{\beta_r} \mu$.

If $r=0$,  then $\la=\mu$ and hence the assertion is trivial.

Let $r>0$. 

\snoi
(a)\ 
Assume that   $\al=\beta_1$.  Then we have $s_\al\la \ble \mu$.  Hence we may assume that $(\al,\mu)\le0$.  Then we have $s_\al \la \ble \mu \ble s_\al\mu$.

\snoi
(b)\ 
Assume that  $\al\neq \beta_1$. Set $\la\rq{}\seteq s_{\beta_1} \la= s_{\beta_2} \cdots s_{\beta_r} \mu \ble \mu$.  
Then $s_\al \la = s_{(s_\al\beta_1)} (s_\al \la\rq{})$.
Since $0< (\beta_1,\la\rq{})= (s_\al\beta_1,s_\al\la\rq{})$
and $s_\al\beta_1\in\prt$, we have $s_\al\la \ble s_\al \la\rq{}$.
If $(\al,\mu)\ge0$, then we have
$s_\al  \la \ble s_\al\la'\ble\mu$, and
if $(\al,\mu)\le0$, then we have
$s_\al\la\ble s_\al\la'\ble s_\al\mu$.
\QED

\begin{lem} \label{lem:wLavLa}
Let $w,v\in \weyl$.
Then the following conditions are equivalent.
\bna
\item$w\ge v$,\label{ita}
\item
$w\La\ble v\La$ for any $\La\in\pwtl$,\label{itb}
\item$w\La_i \ble v\La_i$ for all $i\in I$.\label{itc}
\ee
\end{lem}
\begin{proof}
\eqref{ita}$\Rightarrow$\eqref{itb}\ 
If $w\ge v$, then
there exists a sequence  of positive real roots $\beta_1,\ldots,\beta_r$ 
such that $\bl s_{\beta_{k+1}} s_{\beta_{k+2}} \cdots s_{\beta_r}v\br^{-1}\beta_k\in\prt$
for all $1\le k\le r$ and $w=s_{\beta_1} s_{\beta_2} \cdots s_{\beta_r}v$. 
Then 
\eqn
&&\ang{\beta_k,\;s_{\beta_{k+1}} s_{\beta_{k+2}} \cdots s_{\beta_r}v\La}
=\ang{\bl s_{\beta_{k+1}} s_{\beta_{k+2}} \cdots s_{\beta_r}v\br^{-1}\beta_k,\;\La}\ge0.
\eneqn

\snoi
\eqref{itb}$\Rightarrow$\eqref{itc} is trivial.
Let us show \eqref{itc}$\Rightarrow$\eqref{ita} 
by induction on $\ell(w)$.
If $\ell(w)=0$,  then $w=\id$ so that $\La_i\ble v\La_i$ for all $i\in I$.  Since $v\La_i \ble \La_i$,  we have $v\La_i=\La_i$ for all $i\in I$ so that $v=\id$.

Assume that $\ell(w)>0$.  Take $a\in I$ such that $s_a w < w$.  Then $w^{-1}\al_a \in \Delta_-$ and hence
$(\al_a,w\La_i)\le 0$ for all $i\in I$.  By the assumption and Lemma \ref{lem:lamu},  we have 
either $s_a w\La_i \ble v \La_i$,  $(\al_a,v\La_i) \ge 0$ or
$s_a w\La_i \ble s_a v \La_i$,  $(\al_a,v\La_i) < 0$.

If $s_a v > v$,  then $v^{-1} \al_a \in \Delta_+$ so that $(\al_a,v\La_i) \ge 0$ for all $i\in I$.
Hence $s_aw\La_i \ble v\La_i$ for all $i\in I$.  By induction on $\ell(w)$,  $s_aw\ge v$ so that 
$w> s_a w \ge v$.

If $s_av < v$,  then $v^{-1} \al_a \in \Delta_-$ so that $(\al_a,v\La_i) \le 0$ for all $i\in I$.
Hence  $s_a w\La_i\ble s_a v\La_i$ for all $i\in I$. 
By induction on $\ell(w)$,  $s_aw \ge s_av$ and hence $w\ge v$,  as desired.
\end{proof}

\begin{coro}
Let $i\in I$ and $w,v\in\weyl$.  If $w\ge v$ and $ws_i\La_i \ble vs_i \La_i$,  then
$ws_i\ge vs_i$.
\end{coro}
\begin{proof}
If $j\neq i$,  then $ws_i\La_j =w\La_j$ and $vs_i \La_j =v\La_j$ so that $ws_i\La_j  \ble v s_i\La_j$.
Hence the assertion follows from the lemma above.
\end{proof}

\subsection{Categories $\catC_{w}$ and $ \catC_{w, v}$} \label{Sec: Cwv}
In this subsection, we recall  the categories $\catC_w$, $\catC_{*,v}$ and $\catC_{w,v}$ defined in  \cite{KKOP18}. 

For $M\in \Modg\bl R(\beta)\br$ we define
\eqn
&&\ba{l}
\gW(M) \seteq  \{  \gamma \in  \rlQ_+ \cap (\beta - \rlQ_+)  \mid  e(\gamma, \beta-\gamma) M \ne 0  \}, \\
\gW^*(M) \seteq  \{  \gamma \in  \rlQ_+ \cap (\beta - \rlQ_+)  \mid  e(\beta-\gamma, \gamma) M \ne 0  \}.
\ea
\eneqn

For $w,v \in\weyl$,  we define the full monoidal subcategories of $R\gmod$
by
\eq
&&\ba{l}
\catC_w\seteq\st{M\in R\gmod\Mid \W\,(M)\subset\prtl\cap\, w\nrtl} ,\\
\catC_{*,v}\seteq\st{M\in R\gmod\Mid \W^*(M)\subset\prtl\cap\, v\prtl},\\
\catC_{w,v} \seteq\catC_w \cap \catC_{*,v}.
\ea\label{def:Cw}
\eneq

An ordered pair $(M,N)$ of $R$-modules is called {\em unmixed}
if $$\sgW(M)\cap\gW(N)\subset\{0\}.$$

Assume that $\lambda, \mu \in \weyl \Lambda$ for some $\La \in \pwtl$ and $\lambda \ble \mu$.
Then there exists an object  $\dM(\lambda, \mu)$ in $R(\mu-\la)\gmod$,  called  the \emph{determinantial module}.  (See \cite[Section 3.3]{KKOP21} for the precise definition and more properties of them.)
Note that  $\dM(\la,\mu)$ is an \afr (\cite[Theorem 3.26]{KKOP21}). 
For $\La \in \pwtl$ and $w,v\in \weyl$ with $v \le w$ we have
\eqn
\dM(w\La,\La) \in \catC_w, \quad \dM(w\La,v\La) \in \catC_{w,v}.
\eneqn

\subsection{Localizations of $\catC_{w}$ and $\catC_{w,v}$ via left braiders}
In this subsection we recall  the localizations of the categories $\catC_{w}$, $\catC_{w,v}$ via left braiders studied in \cite{KKOP21, KKOP22}.

Let $L(i)$ denote the one-dimensional  graded self-dual simple module of $R(\al_i)$.
For any simple module $M\in R\gmod$,  there exists a graded left braider $(M, \coR_M,\dphi_M)$ in $R\gmod$ which is \emph{non-degenerate},  that is,  $\coR_M(L(i))\neq 0$  for all $i\in I$  (\cite[Proposition 4.1]{KKOP21}).    
Such a non-degenerate braider is unique up to a constant multiple (\cite[Lemma 4.3]{KKOP21}).

Let $w \in \weyl$ with $I_w=I$,  where 
$ I_w \seteq \set{i\in I}{w\La_i\neq \La_i}$.
The family of graded left braiders in $R\gmod$
\eqn
\st{\bl\dM(w\La_i,\La_i),\coR_{\dM(w\La_i,\La_i)}, \dphi_{\dM(w\La_i,\La_i)}\br}_{i\in I} 
\eneqn
is a real commuting family  (\cite[Proposition 5.1]{KKOP21}).   Moreover it is a family of central objects in the category $\catC_w$.
Note that (see Remark \ref{rem:changes} (ii))
\eq \label{eq:phi_w}
\dphi_{\dM(w\La_i,\La_i)}(\beta)= -(w\La_i+\La_i, \beta) \qt{for any } \ \beta\in \rtl. 
\eneq

Hence there exist localizations of $R\gmod$ and $\catC_w$ via this real commuting  family of graded left braiders  and we denote them by $\bl R\gmod\br[\dM(w\La_i,\La_i)^{\circ -1}; i\in I]$ and 
$\tcatC_w=\catC_w[\dM(w\La_i,\La_i)^{\circ -1}; i\in I]$,  respectively.  
We have a commutative diagram of functors
\eqn
\xymatrix@C=12em{
\catC_w \akew\ar@{>->}[r]_{\iota_w} \ar[d]_{\Phi_w}& R\gmod \ar[d]_{Q_w} \\
\tcatC_w \akew\ar@{>-->}[r]_{\widetilde\iota_w}&\bl R\gmod \br[\dM(w\La_i,\La_i)^{\circ -1}; i\in I]
}
\eneqn
where $\Phi_w$ and $Q_w$ denote the localization functors,   and $\widetilde \iota_w$ is the induced functor from  the inclusion functor $\iota_w$. 

\begin{thm} [{\cite[Theorem 5.9, Theorem 5.11]{KKOP21},  \cite[Theorem 3.9]{KKOP22}}]
\hfill
\bna
\item The functor $\widetilde \iota_w \cl  \tcatC_w \to \bl R\gmod \br[\dM(w\La_i,\La_i)^{\circ -1}; i\in I]$ is an equivalence of categories.
\item The monoidal category $\tcatC_w$ is rigid, that is, every object of $\tcatC_w$ has a left dual and a right dual.
\ee
\end{thm}

Let $v\in \weyl$ such that $v\le w$. 
 In \cite[Section 4]{KKOP22},  it is shown that there exists a real commuting family of  graded left braiders 
\eqn
\st{ \bl \dM(w\La_i,v \La_i),\coR_{\dM(w\La_i,v\La_i)}, \dphi^\ell_{w,v,\La_i}\br }_{i\in I}. 
\eneqn
 in the category $\catC_{*,v}$. 
Moreover it is a family of central objects in the category $\catC_{w,v}$.
If $v=\id$,  then $\coR_{\dM(w\La_i,v\La_i)}$ is nothing but the non-degenerate braider 
$\coR_{\dM(w\La_i,\La_i)}$, so that  this abuse of notation is justified. 
Note that
\eq
\dphi^\ell_{w,v,\La_i}(\beta) = -(w\La_i+v\La_i, \beta) \qt{for any } \ \beta\in \rtl. 
\eneq

Let us denote by $\catC_{*,v}[\dM(w\La_i,v\La_i)^{\circ -1}; i\in I]$ and 
$\tcatC_{w,v}=\catC_{w,v}[\dM(w\La_i,v\La_i)^{\circ -1}; i\in I]$,  the localizations of $\catC_{*,v}$ and $\catC_{w,v}$,  respectively.
Then we have the commutative diagram of functors
\eqn
\xymatrix@C=12em{
\catC_{w,v}\akew \ar@{>->}[r]_{\iota_{w,v}} \ar[d]_{\Phi_{w,v}}& \catC_{*,v} \ar[d]_{\Ql_{w,v}} \\
\tcatC_{w,v}\akew \ar@{>-->}[r]_{\widetilde\iota_{w,v}}&\bl \catC_{*,v}\br[\dM(w\La_i,\La_i)^{\circ -1}; i\in I]
}
\eneqn
where $\Phi_{w,v}$ and $\Ql_{w,v}$ denote the localization functors,   and $\widetilde \iota_w$ is the induced functor from $\iota_{w,v}$ which is the inclusion functor. 
\begin{thm}[{\cite[Theorem 4.5]{KKOP22}}]
The functor $\widetilde \iota_{w,v} \cl  \tcatC_{w,v} \to \bl \catC_{*,v} \br[\dM(w\La_i,\La_i)^{\circ -1}; i\in I]$ is an equivalence of categories.
\end{thm}

\vskip 2em

\section{$\tCwv$ as the right localization of  $\Cw$} \label{Sec: loc rb}
For $\eta,\beta\in \prtl$, denote the functor
$\Res_{\eta,\beta}\cl \Modg\bl R(\eta+\beta)\br\to \Modg\bl R(\eta)\tens R(\beta)\br$ simply by $\Res_{*,\beta}$. Then  for any $V,W \in R\gmod$,  we have 
(\cite[Theorem~2.1]{BKM14})
\eqn
&&V \conv \Res_{*,\beta} (W) \monoto  \Res_{*,\beta} (V \conv W) \qtq
 \Res_{*,\beta} (V \conv W) \epito q^{(\beta,\wt(W))}  \Res_{*,\beta} (V)\conv W.
\eneqn

The following lemma follows from
\cite[Theorem~2.1]{BKM14}. 
\begin{lem} \label{lem:ResComm}
For $X,Y,Z \in \Modg(R)$ and $\beta\in \prtl$, the diagram below is commutative.
$$
\xymatrix{
X\conv  \Res_{*,\beta} (Y\conv Z)\akew \ar@{>->}[r] \ar@{->>}[d]&\Res_{*,\beta} (X \conv Y\conv Z)  \ar@{->>}[d] &
\\
q^{(\beta,\wt(Z))}  X\conv \Res_{*,\beta}(Y)\conv Z\akew  \ar@{>->}[r]&
q^{(\beta,\wt(Z))}  \Res_{*,\beta} (X\conv Y) \conv Z.
}
$$
\end{lem}

\medskip
Let $w \in \weyl$ with $I_w =I$ and $v \le w$.
\begin{prop} \label{prop:right_braider}
For any $\La\in\pwtl$, there exists a morphism in $\catC_w$
 $$\coRr_{\dM(w\La, v\La)}(X) \cl X \conv \dM(w\La,v\La) \to q^{-(\wt(X),w\La+v\La)} \dM(w\La,v\La ) \conv X$$ 
functorial  in $X\in \catC_w$. 
Moreover, if $X$ belongs to $\catC_{w,v}$,  then the morphism $\coRr_{\dM(w\La, v\La)}(X)$ is an isomorphism.
\end{prop}
\begin{proof}
Let $(\Mw,\coR_{\Mw})$ be the non-degenerate left braider 
in $R\gmod$ associated with $\Mw$.
Then we have an isomorphism 
\eqn
q^{A}X \conv \Mw \isoto[\ \coR_{\Mw}(X)^{-1}\ ] \Mw \conv X
\eneqn
functorial in  $X \in \catC_w$,  where $A=(w\La+\La, \wt(X))$ (see \eqref{eq:phi_w}).

Let $\beta = \La-v \La$.
Recall that 
$$\Res_{*,\beta} \Mw \simeq \Mwv\tens \Mv \qt{(including the grading shift).}$$
Let $\al=v\La -w\La$, and $\gamma=-\wt(X)$.
Then we have  a morphism  in $\bl R(\gamma+\alpha)\tens R(\beta)\br\gmod$
\eq 
&&\hs{5ex}q^A(X  \conv \Mwv )\tens \Mv \simeq q^AX \conv \Res_{*,\beta} (\Mw) \nonumber \\
&&\hs{10ex}\monoto q^A\Res_{*,\beta}(X\conv \Mw) 
\isoto[\Res_{*,\beta}(\coR_{\Mw}(X)^{-1})] \Res_{*,\beta}(\Mw \conv X)  \label{eq:rightbraider} \\
&&\hs{10ex}
\epito
q^{-(\beta, \gamma)}
\Res_{*,\beta}(\Mw) \conv X
\simeq 
q^{-(\beta, \gamma)}
(\Mwv\conv X)\tens \Mv. \nonumber 
\eneq 
By applying the functor $\Hom_{R(\beta)\gmod}(\Mv,-)$ we obtain a morphism in $\catC_w$ 
\eqn 
\coRr_{\Mwv}(X)\cl X\conv \Mwv \to  q^{(-\wt(X),w\La+v\La)} \Mwv \conv X
\eneqn
which is functorial in $X\in \catC_w$.

If an $R(\gamma)$-module $X$ belongs to $ \catC_{*,v}$,  then by \cite[Lemma 4.1]{KKOP22},   we have  isomorphisms 
\eq \label{eq:XResM1}
&&X\conv \Res_{*,\beta}(\Mw) \simeq \Res_{*,\beta}(X\conv \Mw)  \qt{and} \\
\label{eq:XResM2} &&\Res_{*,\beta}(\Mw \conv X) \simeq q^{-(\beta,\gamma)}  (\Mwv \conv X)\tens \Mv.
 \eneq
Hence the composition \eqref{eq:rightbraider} is an  isomorphism so that 
the morphism $\coRr_{\Mwv}(X)$ is an isomorphism for any $X\in \catC_{w,v}$,  as desired.
\end{proof}

By \cite[Theorem 4.12]{KKOP18},  for any $\La,\La\rq{}\in \wtl_+$  we have
\eqn 
\La(\dM(w\La\rq{},v\La\rq{}),\Mwv)=(\wt(\dM(w\La\rq{},v\La\rq{})),w\La+v\La)
= (w\La\rq{}-v\La\rq{},w\La+v\La).
\eneqn
The following corollary is a direct consequence of this and Proposition \ref{prop:right_braider}.
\begin{coro}\label{cor:right_baraider}
Let $\phi^r_{w,v,\La_i}(\gamma)\seteq(\gamma,w\La_i+v\La_i)$ for $\gamma \in \rtl$.
Then the family 
$$\st{\bl \dM(w\La_i,v\La_i), \coRr_{\dM(w\La_i,v\La_i),},\phi^r_{w,v,\La_i}\br}_{i\in I}$$ 
is a real commuting family of right graded braiders in the category $\catC_w$.  It is also a family of central objects in $\catC_{w,v}$.
\end{coro}
Note that  
$$\phi^r_{w,v,\La_i} =-\phi^\ell_{w,v,\La_i}.$$

The following theorem gives a characterization of $\catC_{*,v}$.

\begin{thm}\label{thm:C*v}
A simple module $M$ in $R\gmod$ belongs to $\catC_{*,v}$ if and only if 
$$\tLa(M,\Mv)=0 \qt{for all $\La\in \wtl_+$.}$$
\end{thm}
\begin{proof}
If $M\in\catC_{*,v}$,  then $\bl M,\Mv\br$ is unmixed and hence $\tLa(M,\Mv)=0$.

Assume that $\tLa(M,\Mv)=0$ for all $\La \in \wtl_+$.
By \cite[Proposition 1.24]{KKOP18} and \cite[Theorem 2.19]{TW16},  there exist simple modules $X\in \catC_{*,v}$  and $Y\in \catC_{v}$ such that
$M\simeq X\hconv Y$.
Hence  we have 
\eqn 
0=\tLa(M,\Mv) \ge \tLa(Y,\Mv)\ge0
\eneqn
for any $\La\in \wtl_+$,
where  the  first inequality follows from 
\cite[Theorem 2.11 (ii)]{KKOP22}.
On the other hand, we have
\eqn
\tLa(Y,\dM(v\La,\La))
&&=\dfrac{1}{2} \Bigl( \La\bl Y,\dM(v\La,\La) \br  + \bl\wt(Y), \wt(\dM(v\La,\La))\br \Bigr) \\
&&=\dfrac{1}{2} \Bigl( -\La(\dM(v\La,\La),Y )  + \bl\wt(Y), v\La-\La\br \Bigr) \\
&&=\dfrac{1}{2} \Bigl( \bl\wt(Y), v\La+\La \br  + \bl\wt(Y), v\La-\La\br \Bigr) 
=\bl\wt(Y),v\La\br.
\eneqn

It follows that $(v^{-1}\wt(Y), \La)=0$ for all $\La\in \wtl_+$, so that $\wt(Y)=0$. It follows that 
$M\simeq X \hconv \one=X$ and hence  $M$ belongs to $\catC_{*,v}$ as desired. 
\end{proof}
\begin{coro} \label{coro:belongtoCwv1}
A  simple module  $X$  in $\catC_w$ belongs to $\catC_{w,v}$ if and only if
$X$ commutes with $\Mwv$  and
$\La(X,\Mwv)=(\wt(X), w\La+v\La)$ for all $\La\in \wtl_+$.
\end{coro}
\begin{proof}
The ``only if'' part follows from Proposition \ref{prop:right_braider}.
Let us prove  the ``if'' part.

Assume that $X$ commutes with $\Mwv$  and
$\La(X,\Mwv)=(\wt(X), w\La+v\La)$ for all $\La\in \wtl_+$.

Then, we have
\eqn 
\La(X,\Mw) && =-\La(\Mw,X)=(w\La+\La,\wt(X)),
\eneqn
where the first equality comes from the fact that $X\in \catC_w$ so that $X$ commutes with $\Mw$,  and the second comes from  \cite[Corollary 5.10]{KKOP21}.

Since $\Mw=\Mwv\hconv\Mv$ and $X$ commutes with $\Mwv$,  we have 
$\La(X,\Mw) =\La(X,\Mwv)+\La(X,\Mv)$.
Hence we have
\eqn
\La(X,\Mv)&&=
\La(X,\Mw)-\La(X,\Mwv)\\
&&=(w\La+\La,\wt(X))-(w\La+v\La,\wt(X))=
-(v\La-\La,\wt(X)).
\eneqn
It follows that $\tLa(X,\Mv)=0$ and hence $X$ belongs to $\catC_{w,v}$ by Theorem \ref{thm:C*v}.
\end{proof}

\begin{coro} \label{coro:belongtoCwv2}
A  simple module  $X$  in $\catC_{*,v}$ belongs to $\catC_{w,v}$ if and only if
$X$ commutes with $\Mwv$  and
$\La(\Mwv,X)=-(\wt(X), w\La+v\La)$ for all $\La\in \wtl_+$.
\end{coro}
\begin{proof}
Since  the ``only if'' part is obvious,
let us prove the ``if'' part.

Assume that a simple $X$  in $\catC_{*,v}$
commutes with $\Mwv$ and $\La(\Mwv,X)=-(\wt(X), w\La+v\La)$ for all $\La\in \wtl_+$.
Then  we have
\eqn
\La(X,\Mw)&&=\La(X,\Mwv\hconv \Mv)  \\
&&= \La(X,\Mwv)+\La(X,\Mv)\\&&= (\wt(X), w\La+v\La)-(\wt(X),v\La-\La)= (\wt(X),w\La+\La).
\eneqn 
Hence we have
\eqn 
 2 \de(X,\Mw)&&= \La(X,\Mw)+ \La(\Mw,X)\\
&& =(\wt(X),w\La+\La)+ \La(\Mw,X)\\
&&\le (\wt(X),w\La+\La)-(\wt(X),w\La+\La)=0,
\eneqn
where the inequality is \cite[Proposition 4.4]{KKOP21}.
 Hence we have $\La(\Mw,X)=-(\wt(X),w\La+\La)$, and
$\coRl_{\Mw}(X)$ does not vanish by \cite[Proposition 4.4]{KKOP21}.
Since $X$ commutes with $\Mw$, $\coRl_{\Mw}(X)$ is an isomorphism.
Then \cite[Corollary 5.10]{KKOP21} implies that  $X$ belongs to $\catC_w$. 
\end{proof}

\medskip
By Corollary~\ref{cor:right_baraider},  we have  localizations 
\eqn
\catC_w \to \catC_w[\dM(w\La_i,v\La_i)^{\circ -1}; i\in I] \qt{and} \quad \catC_{w,v} \to \catC_{w,v}[\dM(w\La_i,v\La_i)^{\circ -1}; i\in I]. 
\eneqn

Let us denote $\catC_{w,v}[\dM(w\La_i,v\La_i)^{\circ -1}; i\in I]$ by $\tcatC_{w,v}$.
By the definition of localization,  the embedding $\iota_{w,v}\cl \catC_{w,v} \monoto \catC_{w}$ induces a fully faithful functor $$\tilde \iota_{w,v}\cl\tcatC_{w,v} \seteq\catC_{w,v}[\dM(w\La_i,v\La_i)^{\circ -1}; i\in I] \monoto  \catC_{w}[\dM(w\La_i,v\La_i)^{\circ -1}; i\in I].$$
Note that the subcategory $\tcatC_{w,v}$ is closed by taking subquotients and extensions
in $\catC_{w}[\dM(w\La_i,v\La_i)^{\circ -1}; i\in I]$ (\cite[Proposition 2.10]{KKOP21}).

\begin{thm} \label{thm:essential_surj}
The functor $\tilde \iota_{w,v}\cl\tcatC_{w,v} \to  \catC_{w}[\dM(w\La_i,v\La_i)^{\circ -1}; i\in I]$
is an equivalence of monoidal categories.
\end{thm}
\begin{proof} In the course of the  proof, we omit the grading shifts.
Let $Q$ denote the localization functor $\catC_w \to \catC_{w}[\dM(w\La_i,v\La_i)^{\circ -1}; i\in I]$. 
It remains to show that for any $X\in \catC_w$,  the object $Q(X)$ belongs to $\tcatC_{w,v}$.

\snoi(a)\ 
Assume first that $X$  is simple in $\catC_w$ such that $Q(X)\not\simeq 0$ and $X$ commutes with $\Mwv$ for all $\La\in \wtl_+$.
Since $Q(X)\not\simeq 0$,  we have 
$$\coRr_{\Mwv}(X)\neq 0 \qt{for all} \quad \La \in \wtl_+.$$
Since $\Mwv$ is \afr,   we have
$\coRr_{\Mwv}(X)=\rmat{X,\Mwv}$ up to a constant multiple and hence
$$\La(X,\Mwv)=\phi^r_{w,v,\La}(\wt (X))=(\wt(X),w\La+v\La)  \qt{for all} \quad \La \in \wtl_+.$$
Hence $X$ belongs to $\catC_{w,v}$,  by Corollary \ref{coro:belongtoCwv1}.
Thus $Q(X)$ belongs to $\tcatC_{w,v}$.

\snoi(b)\ 
Assume that $X$ is simple in $\catC_w$ such that $Q(X)\not \simeq 0$. 
If $\de(X,\Mwv) >0$ for $\La \in \wtl_+$,  then we have $\de(X\hconv \Mwv,\Mwv) < \de(X,\Mwv) $ by 
\cite[Corollary 3.18]{KKOP21}.  Hence by taking large enough $\la \in \wtl_+$,  we may assume that  $\de(X\hconv \dM(w\la,v\la), \Mwv)=0$ for any $\La\in \wtl_+$.  
Since  $\coRr_{\dM(w\la,v\la)}(X)$ is decomposed into $$X \conv \dM(w\la,v\la) \epito X \hconv \dM(w\la,v\la) \monoto \dM(w\la,v\la) \conv X$$
and $Q(\coRr_{\dM(w\la,v\la)}(X))$ is an isomorphism,
we have  $$Q(X\hconv  \dM(w\la,v\la))\simeq Q(X) \conv  \dM(w\la,v\la). $$
Hence the object
$Q(X)\simeq  Q(X\hconv \dM(w\la,v\la)) \conv  \dM(w\la,v\la)^{\circ -1}$ belongs to $\tcatC_{w,v}$ by (a).

\snoi
(c)\ 
Since the subcategory $\tcatC_{w,v}$ of  $\catC_{w}[\dM(\La_i,\La_i)^{\circ -1}; i\in I]$ is closed under extension,  every object $Q(X)$ for $X$ in $\catC_w$ belongs to $\tcatC_{w,v}$,  as desired.
\end{proof}

Let $\Qr_{w,v}$ denote the composition of functors
\eq
\Qr_{w,v} \cl \catC_w \to \catC_{w}[\dM(w\La_i,v\La_i)^{\circ -1}; i\in I] \isoto \tcatC_{w,v}.\eneq

In the following two propositions,
we characterize the kernels of
$\Ql_{w,v}\cl\catC_{*,v}\to\tcatC_{w,v}$ and $\Qr_{w,v}\cl\catC_{w}\to\tcatC_{w,v}$.

\begin{prop}\label{prop:Laneq1}
Let $X$ be a simple object of $\catC_{*,v}$.
Then, $\Ql_{w,v}(X)\not\simeq 0$ if and only if
$$ \La(\dM(w\la,v\la),X) =-(\wt(X),w\la+v\la)$$
for any $\la \in \wtl_+$.
\end{prop}
\begin{proof}
``Only if'' part is obvious. Let us show the ``if'' part.
 
There exists $\mu \in \wtl_+$ such that $\dM(w\mu,v\mu) \hconv X$ commutes with  $\dM(w\La,v\La)$ for any $\La \in \wtl_+$.
Then we have
\eqn
\La(\dM(w\la,v\la),\dM(w\mu,v\mu) \hconv X)&&=
\La\bl\dM(w\la,v\la),\dM(w\mu,v\mu)\br+\La\bl\dM(w\la,v\la),X\br\\
&&=-\bl w\la+v\la,\wt(\dM(w\mu,v\mu) \hconv X)\br
\eneqn
for any $\la\in\pwtl$.
Hence, Corollary \ref{coro:belongtoCwv2}
implies that
$\dM(w\mu,v\mu) \hconv X \in \catC_{w,v}$.

Then $$ \Ql_{w,v} (\dM(w\mu,v\mu)) \conv \Ql_{w,v}(X)\epito \Ql_{w,v}(\dM(w\mu,v\mu) \hconv X)\not\simeq0, $$
implies that $\Ql_{w,v}(X)\not\simeq0$.
\end{proof}

\begin{prop} \label{prop:Laneq2}
Let $X$ be a simple object of $\catC_{w}$.
Then, $\Qr_{w,v}(X)\not\simeq 0$ if and only if
$$ \La(X,\dM(w\la,v\la))=(\wt(X),w\la+v\la)$$
for any $\la \in \wtl_+$.
\end{prop}
\begin{proof}
The proof is similar to the one of the preceding proposition by using 
Corollary~\ref{coro:belongtoCwv1} instead of Corollary~\ref{coro:belongtoCwv2}.
\end{proof}

\vskip 2em

\section{Properties of $\Cwv$} \label{Sec: properties of Cwv}
\subsection{Right rigidity}
As an application of Theorem~\ref{thm:essential_surj}, 
we will prove the right rigidity of $\tCwv$.
\begin{thm} \label{thm:right_rigidity}
The category $\tcatC_{w,v}$ is right rigid, i.e., every object has a right dual.
\end{thm}
\begin{proof}
In the course of the  proof, we omit the grading shifts.
Let $X \in \catC_{w,v}$.  Since $\catC_{w,v} \subset \catC_w \subset \tcatC_w$ and the category $\tcatC_w$ is right rigid,  there exists $Y\in \catC_w$, $\La\in \wtl_+$ and morphisms in  $\catC_w$
$$X \conv Y \To[\eps] \Mw,\quad \Mw \To[\eta] Y\conv X $$
such that 
the composition 
$$X \conv \Mw \To[X \conv \eta] X  \conv  Y \conv X \To[\eps \conv   X] \Mw\conv X$$
is an isomorphism.

Recall that $\Mw\simeq \Mwv\hconv \Mv$.
Let $\beta=\La-v\La\in\prtl$ and $\gamma= -\wt(X)\in\prtl$.
We have the following commutative diagram:
$$
\xymatrix@C=6ex{
X \conv (\Res_{*,\beta} \Mw) \ar[r] \ar[d]^-{\simeq\; \eqref{eq:XResM1}} & X\conv  \Res_{*,\beta} (Y\conv X) \ar@{>->}[d] & \\
\Res_{*,\beta} (X\conv \Mw) \ar[r]^-{X\circ \eta} & \Res_{*,\beta} (X \conv Y\conv X) \ar[r]^-{\eps \circ X}  \ar@{->>}[d]& \Res_{*,\beta} (\Mw \conv X) \ar@{->>}^{\simeq \;\eqref{eq:XResM2}}[d]\\
& 
\Res_{*,\beta} (X\conv Y) \conv X \ar[r] & 
\Res_{*,\beta} (\Mw)\conv X.
}
$$
Since the composition of the arrows in the middle row is an isomorphism,  we have an isomorphism 
$$X\conv(\Res_{*,\beta}(\Mw)) \isoto\Res_{*,\beta} (\Mw)\conv X. $$
We claim that this isomorphism factors through 
$X\conv \Res_{*,\beta}(Y)\conv X$. Indeed,
in the diagram
\eqn
\xymatrix{
X\conv(\Res_{*,\beta}(\Mw))\ar[r]&X\conv  \Res_{*,\beta} (Y\conv X)\akew \ar@{>->}[r] \ar@{->>}[d]&\Res_{*,\beta} (X \conv Y\conv X)  \ar@{->>}[d] &
\\
&X\conv \Res_{*,\beta}(Y)\conv X \akew  \ar@{>->}[r]&
\Res_{*,\beta} (X\conv Y) \conv X  \ar[d]\\
&&\Res_{*,\beta} (\Mw)\conv X,
}
\eneqn
the square is commutative by Lemma \ref{lem:ResComm}.

Hence we have a sequence of morphisms
\eqn
X\conv(\Res_{*,\beta}(\Mw)) \to 
X\conv \Res_{*,\beta}(Y)\conv X \to 
\Res_{*,\beta} (\Mw)\conv X
\eneqn
whose composition is an isomorphism. 

Applying the functor
$\Hom_{R\gmod}(\Mv,-)$,  we obtain a sequence of morphisms
$$X\conv \Mwv \To[X \conv \eta\rq{}]   
X\conv Y\rq{} \conv X \To[\eps\rq{}\conv X]   
\Mwv\conv X$$
whose composition is an isomorphism
where $Y\rq{}=\Res^{\Mv}(Y)$ and 
$\Res^{\Mv}(-)$ is the functor $\Hom_{R(\beta)\gmod}(\Mv,\Res_{*,\beta} (-))$.
Note that the morphisms $\eps\rq{}$ and $\eta\rq{}$ are given by 
\eqn
\eps\rq{} \cl X\conv Y\rq{} \to\Res^{\Mv}(X\conv Y) \to \Res^{\Mv}(\Mw) \simeq \Mwv
\eneqn
and 
\eqn
\eta\rq{} \cl\Mwv\simeq \Res^{\Mw}(\Mw)  \to \Res^{\Mw}(Y\conv X) \to 
Y\rq{} \conv X. 
\eneqn
Hence the assertion follows by Proposition~\ref{prop:criterion_right_rigidity} below. 
\end{proof}

\Prop \label{prop:criterion_right_rigidity}
Let $\shc$ be an idempotent complete additive monoidal category.  
If there exist morphisms
$$X\tens Y \To[\eps] \one \qt{and} \quad \one \To[\eta] Y\tens X$$
such that the composition 
$$X \To[X\tens \eta] X \tens Y\tens X \To[\eps\tens X] X\quad (\text{respectively, } Y \To[\eta\tens Y] Y \tens X\tens Y \To[Y\tens \eps] Y)$$
is an isomorphism,  then $X$ has a right dual \ro respectively,  
$Y$ has a left dual \/\rf in $\shc$.
\enprop
\begin{proof}
Assume that the composition 
$$f\cl X\To[X\tens \eta] X \tens Y \tens X \To[\eps\tens X] X$$
is an isomorphism. Let us show that $X$ has a right dual.

Let  $g$ be the inverse of $f$
and let $\eta\rq{}$ be the composition
$$\eta\rq{}\cl \one \To[ \eta]  Y \tens X \To[Y\tens g]Y\tens X.$$
Then we have  a commutative diagram
\eqn
\xymatrix{
X \ar[r]^-{X\tens \eta\rq{}} \ar[rd]_-{X\tens \eta} & X\tens Y \tens X  \ar[r]^-{\eps \tens X} &  X \\
& X\tens Y \tens X \ar[r]_-{\eps \tens X} \ar[u]_-{X\tens Y\tens g} & X  \ar[u]_{g}
}
\eneqn
so that the composition of morphisms in the top row is the identity.
 Hence,  by replacing $\eta$ with $\eta\rq{}$,  we may assume 
from the beginning that $f$ is the identity.

Let $p$ be the composition
$$p\cl Y \To[\eta\tens Y] Y\tens X\tens Y \To[Y\tens \eps] Y.$$
Then the following commutative diagram shows that 
$p\circ p=p$:
\eqn
\xymatrix@C=6em{
Y \ar[r]^-{\eta\tens Y} \ar[d]_{\eta\tens Y}\ar@/^2.5pc/[rr]^p& Y\tens X\tens Y  \ar[r]^{Y\tens \eps} \ar[d]_{\eta\tens Y\tens X\tens Y}&  Y \ar[d]_{\eta\tens Y}\ar@/^4pc/[dd]^p \\
Y\tens X \tens Y \ar[r]^-{Y\tens X\tens \eta \tens Y} \ar[rd]_{Y\tens f\tens Y =\id}& Y\tens X \tens Y \tens X \tens Y  \ar[r]^-{Y\tens X \tens Y \tens \eps} \ar[d]_-{Y\tens \eps\tens X\tens Y}& Y\tens X\tens Y \ar[d]_{Y\tens \eps}\\
& Y \tens X\tens Y \ar[r]^{Y\tens \eps} & Y.
}
\eneqn
Let $\tY\seteq\Im p$ so that $p$ is factored as
$Y\epiTo[r] \tY\monoTo[s] Y$ with $r\circ s=\id_{\tY}$. 
Let $\teps$ and $\teta$ be the compositions
\eqn
\teps \cl X \tens \tY \To[X\tens s] X\tens Y \To[\eps] \one \qt{and} \quad
\teta\cl \one \To[\eta] Y\tens X \To[r\tens X] \tY\tens X. 
\eneqn
Then we have the following commutative diagram
\eqn
\xymatrix{
X \ar[rr]^-{X\tens \teta} \ar[dr]_{X\tens \eta} \ar@/_4pc/ [dddrr]_-{X\tens \eta}&& X \tens \tY \tens X \ar[rr]^-{\teps\tens X} \ar[dr]_-{X\tens s \tens X}& &X \\
&X \tens Y \tens X \ar[rr]^-{X\tens p\tens X} \ar[ur]_-{X\tens r \tens X} \ar[dr]_{X\tens \eta\tens Y\tens X} && X \tens Y \tens X \ar[ur]_{\eps\tens X}&&\\
&&X\tens Y\tens X\tens Y\tens X \ar[ur]_{X\tens Y\tens \eps \tens X} && \\
&& X\tens Y \tens X \ar[u]_{X\tens Y \tens X \tens \eta} \ar@/_2pc/ [uur]_-{X\tens Y \tens f} \ar@/_4.5pc/ [uuurr]_{\eps\tens X}  &&
}
\eneqn
so that the  composition in the top row is the identity $\id_X$.

The composition of the middle column of the below commutative diagram
\eqn
\xymatrix@C=4.5em{
Y \ar[rr]^{\eta\tens Y} \ar[d]_{r} & &Y\tens X \tens Y\ar[d]_-{r\tens X\tens r} \ar[dl]_-{Y\tens X\tens r} \ar[dr]^-{r\tens X\tens p}&&  \\
\tY \ar[r]_-{\eta \tens \tY}  & Y\tens X \tens \tY \ar[r]_-{r\tens X\tens \tY} \ar[dr]_{p\tens X\tens s}& \tY \tens X \tens \tY \ar[r]_-{\tY \tens X\tens s }  \ar[d]_{s\tens X\tens s}& \tY\tens X \tens Y \ar[r]_{\tY \tens \eps} \ar[dl]^{s\tens X\tens Y}& \tY \ar[d]_s\\
   & &Y\tens X\tens Y \ar[rr]_{\tY\tens\eps} && Y
}
\eneqn
 is $p\tens X\tens p$ and hence 
we have $$s\circ(\tY\tens \teps)\circ (\teta\tens \tY)\circ r=(Y\tens \eps) \circ (p\tens X\tens p) \circ (\eta\tens Y)=p,$$
where the last equality follows from
the commutative diagram
\eqn
\xymatrix{
Y\ar[dr]_p \ar[rr]_{\eta\tens Y}& &Y\tens X\tens Y \ar[rr]_{p\tens X\tens p}\ar[dr]_{Y\tens X\tens p}&&Y\tens X\tens Y  \ar[rr]_{Y\tens \eps}&& Y \\
& Y \ar[rr]_{\eta\tens Y}  && Y\tens X\tens Y \ar[rr]_{Y\tens \eps} \ar[ur]_{p\tens X\tens Y}&& Y \ar[ur]_p.&
}
\eneqn
Hence we have
$$(\tY\tens \teps)\circ (\teta\tens \tY)=r\circ p\circ s = r\circ s\circ r\circ s =\id_{\tY}, $$
as desired.
\end{proof}

We conjecture that $\tcatC_{w,v}$ is a rigid monoidal category. 

\subsection{Relations among $\tCwv$}

The following is known as \emph{T-systems}. 
Note that the proof  in  \cite[Proposition 3.2]{GLS13}  also works   for
non-symmetric case in a similar way. 
\begin{prop}[{\cite[Proposition 3.2]{GLS13}}]   \label{prop:Tsystem}
Let $i\in I$ and $w,v\in \weyl$ satisfying $w<ws_i$,  and $v<vs_i$. 
Then we have the following equalities in $\Aq[\mathfrak n]$.
\eqn 
&& q^{(w\La_i-v\La_i,ws_i\La_i)}  \mD(w\La_i,v\La_i) \mD(ws_i\La_i,vs_i\La_i)\\
&&\hs{15ex}= q^{\sfd_i+(vs_i\La_i,v\La_i-ws_i\La_i)} \mD(w\La_i,vs_i\La_i) \mD(ws_i\La_i,v\La_i) + \mD(w\la,v\la)\\
&&\hs{15ex}= q^{\sfd_i+(v\La_i,vs_i\La_i-w\La_i)} \mD(ws_i\La_i,v\La_i) \mD(w\La_i,vs_i\La_i) + \mD(w\la,v\la), 
\eneqn
where $\la=s_i\La_i+\La_i$ and $\sfd_i\seteq(\al_i,\al_i)/2$. 
\end{prop}

\begin{prop}  
\label{prop:lamu=Lai}
Let $i\in I$ and $w,v\in \weyl$ satisfying $w\ge v$, $w<ws_i$,  and $v<vs_i$.
\bna
\item
If $w\ge vs_i$,  then
we have a short exact sequence in $R\gmod$ 
\eqn \label{eq:Tsystem}
0\to&& q^{\sfd_i+(vs_i\La_i,v\La_i-ws_i\La_i)}  \dM(w\La_i,vs_i\La_i)\conv \dM(ws_i\La_i,v\La_i) \\
&&\to q^{A} \dM(w\La_i,v\La_i) \conv \dM(ws_i\La_i,vs_i\La_i) \to \dM(w\la,v\la) \to 0
\eneqn
where $A=(v\La_i,vs_i\La_i-ws_i\La_i)=(ws_i\La_i,w\La_i-v\La_i)$
and $\la=s_i\La_i+\La_i.$
\item If $w \not \ge vs_i$,  then $w\La_i \not \ble vs_i\La_i $ and we have
\eqn
q^{A} \dM(w\La_i,v\La_i) \conv \dM(ws_i\La_i,vs_i\La_i) \simeq \dM(w\la,v\la).
\eneqn
\end{enumerate}
Hence in the both cases   we have
 \eqn 
 q^{(ws_i\La_i,w\La_i-v\La_i)} \dM(w\La_i,v\La_i) \hconv \dM(ws_i\La_i, vs_i\La_i)  \simeq \dM(w(\La_i+s_i\La_i),v(\La_i+s_i\La_i)).
\eneqn
\end{prop}
\begin{proof}
Since $w\ge v$,  we have $w\La_j \ble v\La_j = vs_i\La_j$ for all $j\neq i$.
If $w \not \ge vs_i$,   then we get $w\La_i \not \ble vs_i\La_i$  by Lemma \ref{lem:wLavLa}.
Then $\mD(w\La_i, vs_i\La_i)=0$ and hence Proposition \ref{prop:Tsystem} implies (b).

Assume that  $w \ge vs_i$.  Then
$q^{(vs_i\La_i,v\La_i-ws_i\La_i)}  \dM(w\La_i,vs_i\La_i)\conv \dM(ws_i\La_i,v\La_i) $ 
is a simple module and it is self-dual.
 Thus Proposition~\ref{prop:Tsystem}
and \cite[Lemma 3.2.18]{KKKO18} imply (a), as desired. 
\end{proof}

\begin{thm} \label{thm:generalT}
Let $w \ge v$, $ws_i>w$,  $vs_i>v$,  and
 $\la,\mu\in \wtl_+$.  
\bnum
\item
Either $\xi\seteq\la+s_i\mu\in \wtl_+$ or $s_i\xi=s_i\la+\mu\in \wtl_+$.
\item We have
\eqn
q^{(w\la-v\la, ws_i\mu)}\dM(w\la,v\la) \hconv \dM(ws_i\mu,vs_i\mu) \simeq \dM(w\xi, v \xi).
\eneqn
\item We have \eqn
\tLa\bl \dM(w\la,v\la), \dM(ws_i\mu,vs_i\mu)\br&&=(w\la-v\la,ws_i\mu)=
-(v\la,ws_i\mu-vs_i\mu),\\
\La\bl \dM(w\la,v\la), \dM(ws_i\mu,vs_i\mu)\br&&=(w\la-v\la,ws_i\mu+vs_i\mu)=
-(w\la+v\la,ws_i\mu-vs_i\mu).
\eneqn
\item If $w\not \ge vs_i$,  then $\dM(w\la,v\la)$ and $\dM(ws_i\mu,vs_i\mu)$ commute.
\end{enumerate}
\end{thm}
\begin{proof}
(i)\ For $j\in I$,  we have 
$\ang{h_j,\la+s_i\mu} = \ang{h_j,\la+\mu}-\ang{h_i,\mu}\ang{h_j,\al_i}$ so that 
$\ang{h_j,\la+s_i\mu} \ge 0$ for $j\neq i$.
Since $\ang{h_i,\la+s_i\mu} = \ang{h_i,\la}-\ang{h_i,\mu}$ and 
$\ang{h_i,s_i\la+\mu} = \ang{h_i,\mu}-\ang{h_i,\la}$ we have either $\xi\in \wtl_+$ or $s_i\xi \in \wtl_+$.
\smallskip

(ii)\ In the proof,  we omit the grading shifts.
Set $\dC_\la\seteq\dM(w\la,v\la)$ and $\dC\rq{}_\mu\seteq\dM(ws_i\mu,vs_i\mu)$.
It is enough to show that there is an epimorphism 
$\dC_\la \conv  \dC\rq{}_\mu \epito \dM(w\xi,v\xi)$.

(1) Assume that  $\la=a \La_i$ and $\mu=b\La_i$ for some $a,b\in \Z_{\ge 0}$.  We may assume that $a,b>0$.
We will proceed by induction on $a+b$.  Set $\la\rq{}\seteq(a-1)\La_i$ and $\mu\rq{}\seteq(b-1)\La_i$. 
Note that  $\eta\seteq\La_i+s_i  \La_i\in \wtl_+$ and $\dC_\eta= \dC\rq{}_\eta$.
Hence we have
\eqn
&&\dC_\la \conv \dC\rq{}_\mu \simeq \dC_{\la\rq{}} \conv \dC_{\La_i} \conv \dC\rq{}_{\La_i} \conv \dC\rq{}_{\mu\rq{}}
\epito 
 \dC_{\la\rq{}} \conv \dC_\eta  \conv \dC\rq{}_{\mu\rq{}}\\
&& \hskip 4em\simeq 
 \dC_\eta \conv  \dC_{\la\rq{}}  \conv \dC\rq{}_{\mu\rq{}} 
\epito 
 \dC_\eta \conv  \dC_{\la\rq{} +\mu\rq{}}
\epito 
\dC_{\eta+\la\rq{} +\mu\rq{}},
\eneqn
where  the first epimorphism follows from Proposition~\ref{prop:lamu=Lai} and 
the second last epimorphism follows from the induction hypothesis.

(2) Set $\la=\la\rq{}+a\La_i$, $\mu=\mu\rq{}+b\La_i$,  and $\eta\rq{}\seteq a\La_i+bs_i\La_i$,  where $a=\ang{h_i,\la}$ and $b=\ang{h_i,\mu}$.
Then we have
\eqn 
\dC_\la \conv \dC\rq{}_\mu \simeq
 \dC_{\la\rq{}} \conv \dC_{a\La_i} \conv\dC\rq{}_{b\La_i} \conv \dC\rq{}_{\mu\rq{}}
\epito 
 \dC_{\la\rq{}} \conv \dC_{\eta\rq{}} \conv \dC\rq{}_{\mu\rq{}}.
\eneqn
Since $ \dC_{\la\rq{}}= \dC\rq{}_{\la\rq{}}$ and $\dC\rq{}_{\mu\rq{}}=\dC_{\mu\rq{}}$, 
we have 
\eqn
\dC_{\la\rq{}} \conv \dC_{\eta\rq{}} \conv \dC\rq{}_{\mu\rq{}}\simeq
\begin{cases}
\dC_{\la\rq{}+\eta\rq{}+\mu\rq{}} \simeq \dC_{\la+s_i\mu} = \dM(w\xi,v\xi)& \qt{if} \quad \eta\rq{} \in \wtl_+,\\
\dC\rq{}_{\la\rq{}+\eta\rq{}+\mu\rq{}}\simeq \dC\rq{}_{s_i\la+\mu} = \dM(w\xi,v\xi)& \qt{if} \quad s_i\eta\rq{} \in \wtl_+,
\end{cases}
\eneqn
as desired.

\snoi
(iii) follows from (ii) and \cite[Lemma 3.1.4]{KKKO18}.

\snoi
(iv)\ Since $\dM(w\la,v\la)$ is a product of $\dM(w\La_j,v\La_j)$\rq{}s, and $\dM(ws_i\mu,vs_i\mu)$ is a product of $\dM(w\La_k,v\La_k)$ ($k\not=i$) together with $\dM(ws_i\La_i,vs_i\La_i)$,  the assertion follows from Corollary \ref{prop:lamu=Lai} (b).
\end{proof}
Recall the functors
\eqn
&&\Ql_{w,v}\cl\catC_{*,v}\to\tcatC_{w,v}\qtq\Qr_{w,v}\cl\catC_{w}\to\tcatC_{w,v}.
\eneqn
\Cor\label{cor:Qcenter}
 Let $w \ge v$, $ws_i>w$,  $vs_i>v$,  and
 $\la,\mu\in \wtl_+$. 
\bnum
\item If $\xi\seteq \la+s_i\mu\in\pwtl$, then we have
$$\Ql_{w,v}\bl\dM(ws_i\mu,vs_i\mu)\br\simeq 
q^{-(w\la-v\la,ws_i\mu)}\dM(w\la,v\la)^{\circ-1}
\conv \dM(w\xi,v\xi).$$
\item If $\eta\seteq s_i\la+\mu\in\pwtl$, then we have
$$\Qr_{ws_i,vs_i}\bl\dM(w\la,v\la)\br\simeq 
q^{-(w\la-v\la,ws_i\mu)}\dM(ws_i\eta\,vs_i\eta)
\conv \dM(ws_i\mu,vs_i\mu)^{\circ-1}.$$
\ee
\encor
\Proof
Since the proof is similar, we only prove (i).

By Theorem \ref{thm:generalT} (iii)
and Proposition~\ref{prop:Laneq1}, we have
$\Ql_{w,v}\bl\dM(ws_i\mu,vs_i\mu)\br\not\simeq0$.
Hence 
$$\coRl_{\dM(w\la,v\la)}\bl\dM(ws_i\mu,vs_i\mu)\br
\cl \dM(w\la,v\la)\conv\dM(ws_i\mu,vs_i\mu)\To\dM(ws_i\mu,vs_i\mu)\conv\dM(w\la,v\la)$$
does not vanish.
Since it is an isomorphism in $\tcatC_{w,v}$,
its image $\dM(w\la,v\la)\hconv\dM(ws_i\mu,vs_i\mu)
\simeq \dM(w\xi,v\xi)$ (in $R\gmod$)
is isomorphic to $\dM(w\la,v\la)\conv\Ql_{w,v}\bl\dM(ws_i\mu,vs_i\mu)\br$ in  $\tcatC_{w,v}$.
\QED

\begin{thm} \label{thm:Cwv=Cwsivsi}
Let $i\in I$ and $w,v\in \weyl$ satisfying $v< w$, $w<ws_i$,  and $v<vs_i$.
If $w\not \ge vs_i$,  then we have $\catC_{w,v}=\catC_{ws_i,vs_i}$.
\end{thm}
\begin{proof}
Set $\la=s_i\La_i+\La_i\in\pwtl$.
 Note that   $\Mwv$ and $\dM(w s_i\La_i, vs_i\La_i)$ commute 
and $\Mwv\conv\dM(w s_i\La_i, vs_i\La_i)\simeq\dM(w \la, v\la)$ by Proposition \ref{prop:lamu=Lai} (b).

Assume that a simple module $X$ belongs to $\catC_{w,v}$.  
Then  we have
$X\in \catC_{ws_i}$ and hence 
\eqn
(\wt(X), w\la+v\la)=\La(X,\dM(w\la,v\la)) 
&&= \La(X,\dM(w\La_i,v\La_i))+\La(X,\dM(ws_i\La_i,vs_i\La_i))\\
&&=(\wt(X), w\La_i+v\La_i)+\La(X,\dM(ws_i\La_i,vs_i\La_i))
\eneqn
so that 
$$\La(X,\dM(ws_i\La_i,vs_i\La_i)) = (\wt(X), w(\la-\La_i)+v(\la-\La_i))
=(\wt(X), ws_i\La_i+vs_i\La_i).$$
Because $X\in \catC_{w,v}$,  we have 
\eqn
0=\de\bl X,\dM(w\la,v\la)\br=\de\bl X, \dM(w\La_i,v\La_i)\br+\de\bl X,\dM(ws_i\La_i,vs_i\La_i)\br=\de\bl X,\dM(ws_i\La_i,vs_i\La_i)\br.
\eneqn
Hence by Corollary \ref{coro:belongtoCwv1},  $X$ belongs to $\catC_{ws_i,vs_i}$.

If $X$ is a simple module  in $\catC_{ws_i,vs_i}$,  then  $X$ belongs to $\catC_{w,v}$ by the same argument as the above using Corollary \ref{coro:belongtoCwv2}.

Since the categories $\catC_{w,v}$ and $\catC_{ws_i,vs_i}$ are closed under extensions,  we obtain that $\catC_{w,v}=\catC_{ws_i,vs_i}$,  as desired.
\end{proof}

\begin{thm} \label{thm: equiv between tcwv}
If $w \ge v$, $ws_i > w$  and $vs_i >v$,  then
there is an equivalence of monoidal categories
$$\tcatC_{w,v} \simeq \tcatC_{ws_i,vs_i}.$$
\end{thm}

\begin{proof}
Set $\shc=\catC_{ws_i,v}$.
Then, $\catC_{w,v}\subset\shc$ and  $\catC_{ws_i,vs_i}\subset\shc$.
In  $\shc$,  there exist
a real commuting family of left braiders
$\st{\dM(w\La_j,v\La_j), \coRl_{\dM(w\La_j,v\La_j),},\phi^l_{w,v,\La_j} }_{j\in I}$
and a   real commuting family of right braiders
$\st{\dM(ws_i\La_j,vs_i\La_j), \coRr_{\dM(ws_i\La_j,vs_i\La_j),},\phi^r_{ws_i,vs_i,\La_j} }_{j\in I}$.

Let us denote by
\eqn
\shc^\ml&&\seteq\catC_{ws_i,v}[\dM(w\La_j,v\La_j)^{\circ-1};
j\in I],\\
\shc^\mr&&\seteq\catC_{ws_i,v}[\dM(ws_i\La_j,vs_i\La_j)^{\circ-1};
j\in I]
\eneqn
their localizations, and denote the corresponding localization functors by $\LQ: \shc \to \shc^\ml$ and $\RQ: \shc \to \shc^\mr$,  respectively. 
Since the composition of the fully faithful functors
$$\tcatC_{w,v}\To \shc^\ml\To\catC_{*,v}[\dM(w\La_j,v\La_j)^{\circ-1};
j\in I]$$
is an equivalence, 
the functor $\tcatC_{w,v}\To \shc^\ml$ is an equivalence of monoidal categories.
Similarly,
since the composition of the fully faithful functors
$$\tcatC_{ws_i,vs_i} \To \shc^\mr\To\catC_{ws_i}[\dM(ws_i\La_j,vs_i\La_j)^{\circ-1};
j\in I]$$ 
is an equivalence,  the functor 
$\tcatC_{ws_i,vs_i}\To \shc^\mr$ is an equivalence of monoidal categories.
Hence it is enough to show that
$\shc^\ml$ and $\shc^\mr$ are equivalent as monoidal categories.

In order to see this, we shall prove that
\eq
\LQ\qt{factors as}
\xymatrix@C=7ex{\shc\ar[r]_{\RQ}\ar@/^1.5pc/[rr]|{\LQ}&\shc^\mr\ar@{.>}[r]_{\Phi}&\shc^\ml},\label{phipsi1}\\
\RQ\qt{factors as}
\xymatrix@C=7ex{\shc\ar[r]_{\LQ}\ar@/^1.5pc/[rr]|{\RQ}&\shc^\ml\ar@{.>}[r]_{\Psi}&\shc^\mr}.\label{phipsi2}
\eneq
Since the proof of \eqref{phipsi2} is similar, we shall prove only
\eqref{phipsi1}.
Set $\dC'_\mu=\dM(ws_i \mu,vs_i\mu)$ for any $\mu\in\pwtl$.
By Theorem~\ref{Thm: graded localization}, it is enough to show that
\bna
\item $\LQ(\dC'_\mu)$ is invertible in $\shc^\ml$ for any $\mu\in\pwtl$,
\item $\LQ\bl M\conv\dC'_\mu\br\To[\LQ(\coRr_{\dC'_\mu})]\LQ\bl\dC'_\mu\conv M\br$
is an isomorphism in $\shc^\ml$ for any $\mu\in\pwtl$ and $M\in\shc$.
\ee

\smallskip
(a) follows from Corollary~\ref{cor:Qcenter}.

\smallskip
Let us show (b).
Let $0\to Z\to M\conv\dC'_\mu\to \dC'_\mu\conv M\to Z'\to0$
be an exact sequence.
Since $\RQ\bl\coRr_{\dC'_\mu}(M)\br$ is an isomorphism, we have
$\RQ(Z)\simeq\RQ(Z')\simeq0$.

Then Lemma~\ref{lem:QrQl} below implies that
$\LQ(Z)\simeq\LQ(Z')\simeq0$ and hence
$\LQ\bl\coRr_{\dC'_\mu}(M)\br$ is an isomorphism.

\smallskip
Thus there exist functors
$\Phi\cl \shc^\mr\To\shc^\ml$ and
$\Psi\cl \shc^\ml\To\shc^\ml$, and it is obvious that they are quasi-inverse to each other.
\QED

\begin{lem} \label{lem:QrQl}
Assume that $w \ge v$, $ws_i > w$,  and $vs_i >v$. 
Let  $Z \in \catC_{ws_i,v}$.
Then 
 $\Qr_{ws_i,vs_i}(Z)\simeq 0$ if and only if $\Ql_{w,v}(Z)\simeq 0$. \end{lem}
\begin{proof}
Since the proof is similar,   we only  prove  ''only if'' part. 
We may assume that $Z$ is simple.

Assuming that $\Qr_{ws_i,vs_i}(Z)\simeq 0$ and $\Ql_{w,v}(Z)\not\simeq0$,
we shall derive a contradiction.

\smallskip
Let us denote $\dC_\la \seteq\dM(w\la,v\la)$ and $\dC\rq{}_\mu\seteq\dM(ws_i\mu,vs_i\mu)$ for $\la,\mu\in \wtl_+$.
Then, 
$\coRl_{\dC_\la}(Z)$ does not vanish for any $\la\in\pwtl$.
Hence $\dC_\la\hconv Z\simeq\Im\bl\coRl_{\dC_\la}(Z)\br$
is isomorphic to $\dC_\la\conv Z$ is $\tcatC_{w,v}$, and hence
$\Ql_{w,v}(\dC_\la\hconv Z)\not\simeq0$.

There exists $\la_0\in\pwtl$ such that
$Z'\seteq\dC_{\la_0}\hconv Z$ commutes with $\dC_\La$ for any $\La\in\pwtl$.
We have $\Qr_{ws_i,vs_i}(Z')\simeq 0$ and $\Ql_{w,v}(Z')\not\simeq0$.
By replacing $Z$ with $Z'$,  we may assume from the beginning that $Z$ commutes with all $\dC_\La$ for $\La\in \wtl_+$.  
By Proposition~\ref{prop:Laneq1} and $\Ql_{w,v}(Z)\not\simeq0$,
we have
$\La(\dC_\la,Z)=-(w\la+v\la,\wt(Z))$ for any $\la\in\pwtl$.
Hence Corollary~\ref{coro:belongtoCwv2} implies that
$Z\in\catC_{w,v}$.
Hence we have 
$$\text{$\La(Z,\dC_\la)=(w\la+v\la,\wt(Z))$ for any $\la\in\pwtl$.}$$

Since $\Qr_{ws_i,vs_i}(Z)\simeq 0$, Proposition~\ref{prop:Laneq2}
implies that there exists $\mu\in\pwtl$ such that
$$\La(Z,\dC'_\mu)\not=\bl ws_i\mu+vs_i\mu,\wt(Z)\br.$$

Let us take $\la\in \wtl_+$ such that $\xi\seteq\la+s_i\mu\in\pwtl$.
Then we have a contradiction
\eqn
(w\xi+v\xi ,\wt(Z))
&&=\La(Z,\dC_\xi)= \La(Z,\dC_\la \hconv \dC'_\mu)  \\
&&= \La(Z,\dC_\la) +\La(Z, \dC'_\mu) 
=(w\la+v\la, \wt(Z))+\La(Z, \dC'_\mu) \\
&&\neq (w\la+v\la, \wt(Z))+(ws_i\mu + vs_i\mu, \wt(Z) ) \\
&&=(w\xi+v\xi, \wt(Z)).
\eneqn
Here
the third equality follows from the commutativity of $Z$ and $\dC_\la$.
\end{proof}

\end{document}